\documentclass[11pt]{article}

\usepackage{graphicx}
\usepackage{latexsym,amsmath,amsfonts,amscd, amsthm, dsfont}
\usepackage{bm,color}
\usepackage{epsfig,verbatim,epstopdf,graphics}
\usepackage{subfigure}
\usepackage{changebar}
\usepackage{multirow}

\usepackage{algorithmic}

\usepackage{url}    

\usepackage{yhmath}
\usepackage{booktabs} 
\usepackage{tikz}
\usepackage{verbatim}
\usetikzlibrary{arrows,backgrounds,snakes,shapes}
\numberwithin{equation}{section}

\usepackage{xcolor}   

\graphicspath{{./}{./figure/}}
\allowdisplaybreaks

\newtheoremstyle{plainNoItalics}{}{}{\normalfont}{}{\bfseries}{.}{ }{}

\theoremstyle{plain}
\newtheorem{thm}{Theorem}[section]

\theoremstyle{plainNoItalics}

\newcommand{\mD}{{\mathcal D}}

\newcommand{\be}{\begin{eqnarray}}
\newcommand{\ee}{\end{eqnarray}}
\newcommand{\beno}{\begin{eqnarray*}}
	\newcommand{\eeno}{\end{eqnarray*}}

\makeatletter

\newcommand{\Rmnum}[1]{\expandafter\@slowromancap\romannumeral #1@}
\makeatother

\newcommand{\nv}{\mathbf{ n}}
\newcommand{\uv}{\mathbf{ u}}
\newcommand{\qv}{\mathbf{ q}}
\newcommand{\Av}{\mathbf{ A}}
\newcommand{\Bv}{\mathbf{ B}}
\newcommand{\Fv}{\mathbf{ F}}

\newcommand{\Ec}{{\mathcal E}}
\newcommand{\ptot}{p_{\text{tot}}}

\begin{document}
\baselineskip=1.8pc


\begin{center}
{\bf
 A Kernel Based High Order \lq\lq Explicit'' Unconditionally Stable Constrained Transport Method for Ideal Magnetohydrodynamics 
}
\end{center}

\vspace{.2in}
\centerline{
Andrew Christlieb\footnote{
 Department of Computational Mathematics, Science and Engineering, Department of Mathematics  and  Department of Electrical Engineering, Michigan State University, East Lansing, MI, 48824. E-mail: christli@msu.edu
},
Firat Cakir \footnote{Department of Mathematics, Michigan State University, East Lansing, MI, 48824. E-mail: cakirfir@msu.edu},
and 
Yan Jiang\footnote{School of Mathematical Sciences, University of Science and Technology of China, Hefei, Anhui 230026, People's Republic of China. E-mail: jiangy@ustc.edu.cn }
}

\bigskip
\noindent
{\bf Abstract.}

The ideal Magnetohydrodynamics (MHD) equations are challenging because one needs to maintain the divergence free condition, $\nabla \cdot \Bv = 0$. Many numerical methods have been developed to enforce this condition. In this work, we further our work on mesh aligned constrained transport by developing a new kernel based approach for the vector potential in 2D and 3D. The approach for solving the vector potential is based on the method of lines transpose and is A-stable, eliminating the need for diffusion limiters needed in our previous work in 3D. The work presented here is an improvement over the previous method in the context of problems with strong shocks due to the fact that we could eliminate the diffusion limiter that was needed in our previous version of constrained transport. The method is robust and has been tested on the 2D and 3D cloud shock, blast wave and field loop problems.

%


\vfill


\newpage

\section{Introduction}

The ideal Magnetohydrodynamics  (MHD) equations are one of the most important classical models of plasma physics explaining the macroscopic phenomena of a quasi-neutral plasma system. The ideal MHD equations are the set of transport evolution equations for the quantities of mass, momentum, and energy density as well as the magnetic field in a conducting fluid. Mathematically, the MHD equations are a system of nonlinear hyperbolic conservation laws with divergence-free magnetic field condition. Satisfying this condition in numerical simulations is the main challenge for the numerical methods. Many standard schemes fail to guarantee $\nabla \cdot \Bv = 0$. There have been overall four different dominant approaches overcoming the difficult in the literature: 8-wave formulation \cite{powell1994approximate, powell1999solution}; projection methods \cite{balsara2004comparison, toth2000b}; hyperbolic divergence cleaning methods \cite{dedner2002hyperbolic}; and constrained transport methods \cite{balsara2004second, balsara1999staggered, dai1998simple, evans1988simulation, fey2003constrained, helzel2011unstaggered, londrillo2000high, londrillo2004divergence, rossmanith2006unstaggered, ryu1998divergence, de2001multi, toth2000b, torrilhon2005locally, christlieb2014finite}.

\par The first methodology used for divergence free condition is developed in \cite{brackbill1980effect} utilizing the classical projection method. They solve a Poisson equation to project the incorrect magnetic field to a divergence free subspace using the Hodge decomposition. However, it is difficult  to extend the method to an Adaptive Mesh Refinement (AMR) approach, since this method has to solve a Poisson equation on each time step. The AMR idea is to refine the mesh based on an accuracy criterion, providing increased resolution at shock and increased computational efficiency. Zachary \cite{zachary1994higher} presented a Riemann solver which removes negative pressures and densities and used the Projection method to get a divergence free magnetic field. The 8-wave scheme is the second approach, developed by Powell in \cite{powell1994approximate}. He adds an extra source term to update the magnetic field which satisfies the divergence free condition. By this treatment, the ideal MHD equations, with constraint, becomes a $8\times8$ hyperbolic system with a source term. This scheme is robust and can easily be extended to an AMR framework. However, it has been shown in \cite{toth2000b} that this method can create inaccurate jumps in a discontinuous example, such as rotated shock tube problem, since this procedure is non-conservative because of the source term. The third method introduced by Dedner \cite{dedner2002hyperbolic} is the hyperbolic divergence-cleaning method. This method is a similar to the projection method. In this method, hyperbolic and parabolic corrections are combined together and solved for the magnetic field divergence error. This method is fully explicit, efficient and fast. However, this method has two tunable parameters, the speed of propagation of the error and the damped divergence error rate, which needed to be adjusted to ensure good solutions. This method gives a damped hyperbolic equation for the divergence error and it does not exactly satisfy divergence free condition.
\par The main focus in this work is the constrained transport (CT) method for correcting the magnetic field to satisfy the divergence free constraint. This methodology was invented by Evans and Hawley in \cite{evans1988simulation}. The original CT method is considered to be a modification of the Yee method \cite{yee1966numerical} from electromagnetics for the ideal MHD equations. In the original constrained transport methodology,  staggered electric and magnetic fields are used to create appropriate finite difference operators. These operators eventually lead to a globally divergence-free magnetic field.
\par In the literature, various modification of the constrained transport method have been presented. In particular, high resolution shock capturing schemes have been a primary focus.
 DeVore \cite{devore1991flux} presented an application of a flux corrected transport approach satisfying a divergence free magnetic field. There are a range of different approaches for building the electric field using Ohm's law in a constrained transport methodology including those presented  by Balsara and Spicer \cite{balsara1999staggered}, Dai and Woodard \cite{dai1998simple}, and Ryu et al. \cite{ryu1998divergence}. Londrillo and Zanna \cite{londrillo2000high, londrillo2004divergence} constructed one of the first high order upwind scheme based on the work of Evans and Hawley. De Sterck \cite{de2001multi} introduced a similar constrained transport scheme on unstructured triangle grids based on multidimensional upwind advection schemes. Balsara \cite{balsara2001divergence} described a divergence free adaptive mesh refinement (AMR) method utilizing a constrained transport approach. T{\'o}th \cite{toth2000b} compared several of these schemes maintaining divergence free condition and also showed a staggered magnetic field is unnecessary as well as developing some unstaggered CT frameworks.
\par  Unstaggered Constrained Transport schemes are getting increased attention over the last few years, since it is easy to implement and to extend to adaptive mesh refinement (AMR) methods. For example , Fey and Torrilhon \cite{fey2003constrained} developed an unstaggered upwind method satisfying the divergence free constraint. Rossmanith \cite{rossmanith2006unstaggered} designed an ustaggered wave propagation scheme for MHD flows based on the algorithms in \cite{leveque1997wave} using a constrained transport method to keep the divergence free magnetic field. Helzel et al. \cite{helzel2011unstaggered, helzel2013high}  generalized the 2D unstaggered CT work to 3D MHD equations so that the method is applicable on both Cartesian and rectangular mapped grids.
\par In resent years, several high order methods have been developed for the ideal MHD equations using various discretization approaches. Balsara \cite{balsara2009divergence} designed a third order divergence-free weighted essentially non-oscillatory (WENO) methods for MHD equations with third order Runge-Kutta time integration using a staggered magnetic field. Balsara et al. \cite{balsara2013efficient, balsara2009efficient} introduced a high accuracy ADER-WENO schemes for divergence free magnetohydrodynamics on structured meshes, again using a staggered magnetic field. Li et al. \cite{li2011central, fu2018globally} and Cheng et al. \cite{cheng2013positivity} presented high order central discontinuous Galerkin schemes satisfying the divergence free constraint globally for ideal MHD simulations on two overlapping meshes, called primal and dual meshes, by utilizing different discretization for magnetic induction equations. Kawai \cite{kawai2013divergence} introduced a divergence-free high order accurate finite difference scheme which has an effective shock capturing capability for the MHD equations by constructing artificial diffusion terms to capture numerical discontinuous in the magnetic field. 
\par In our previous paper \cite{christlieb2014finite}, we introduced a high order FD-WENO for the ideal MHD in 2D and 3D by developing a high order unstaggered constrained transport methodology for Hamilton Jacobi equations using a version of FD-WENO. In that work, WENO formulation is applied to the central derivative of the solution $\Av$ instead of the flux values on grid points to approximate the one-sided partial derivative terms $A_x^-$ and $A_x^+$ appear in HJ equation. If explicit time stepping is used for HJ equation, then it looks like a convection reaction equation and it tends to be unstable. That's why we needed to add artificial resistivity terms for the 3D case to stabilize it and to control unphysical oscillations in the magnetic field.

\par In this current work, we further our work \cite{christlieb2014finite} on mesh aligned constrained transport by developing  a new kernel-based approach for the magnetic vector potential to solve the ideal MHD equations in 2D and 3D. 
The approach is based on a kernel-based numerical scheme \cite{christlieb2017kernel, christlieb2019kernel},
which is derived from the Method of Lines Transpose (MOL$^T$) \cite{causley2013method, causley2014method, causley2014higher, causley2016method, christlieb2016weno, causley2017method}. The equation is first discretized in time with an explicit method and transformed to a boundary value problem (BVP) at discrete time levels. The spatial operators are converted from local representations, i.e.. derivatives, to global representations using convolutions with kernels, i.e. Green's functions, and is similar in spirt of taking a Fast Fourier transform. Thus our method is fully implicit since the partial derivative terms are represented by global operators using convolution integrals which are communicating all the previous and future data. By this methodology, we update the predicted magnetic field obtained from the base scheme by a corrected divergence free magnetic field. The approach for solving the vector potential is derived from the method of lines transpose and is A-stable, eliminating the need of diffusion limiters introduced in our previous work in 3D \cite{christlieb2014finite} since the kernel based method is a fully implicit method. The corrected magnetic field is computed by 4th-order accurate central finite difference operators that approximates the curl of the magnetic vector potential. The magnetic vector potential is made to satisfy a weakly hyperbolic system using the Weyl gauge condition. This system is solved using our kernel based scheme developed for Hamilton-Jacobi equations. Our solver is coupled with the 5th-order FD-WENO scheme of Jiang and Shu \cite{jiang1996efficient} as the base scheme for ideal MHD equations. Third order explicit strong-stability-preserving (SSP) Runge-Kutta (RK) is used for the time discretization. The work presented here is an improvement over the previous method in the context of problems with strong shocks due to the fact that, using the kernel-based approach, we could eliminate the diffusion limiter that was needed in our previous version of constrained transport. This method is robust and has been tested on the 2D and 3D cloud shock, blast wave and field loop problems.

\par The rest of this paper is organized as follows: in section 2, we briefly review the MHD equations and the evolution of the magnetic vector equations; in section 3, we give a brief outline of the CT algorithm; in section 4, we present our novel numerical scheme for 1D Hamilton-Jacobi equations; we introduce the multidimensional solver in section 5; the resulting 2D and 3D schemes are tested on several numerical problems in section 6. 

\section{The Ideal MHD Equations} 

In this section we present a brief review of the ideal MHD equations, which is a first order hyperbolic system of conservation laws. The conservative form of the ideal MHD equations can be written as 
\begin{equation}
\label{eq:MHD}
\begin{aligned}
	 {\partial_t}
	\begin{bmatrix}
	\rho \\
	\rho \uv \\
	\Ec \\
	\Bv
	\end{bmatrix}
	+ \nabla \cdot
	\begin{bmatrix}
	\rho \uv \\
	\rho \uv \otimes \uv + \ptot \mathbb{I} - \Bv \otimes \Bv \\
	\uv(\Ec + \ptot) - \Bv(\uv \cdot \Bv) \\
	\uv \otimes \Bv - \Bv \otimes \uv
	\end{bmatrix}
	= 0, 
\end{aligned}
\end{equation}
\begin{align}
\label{eq:DivergenceFree}
	\nabla \cdot \Bv = 0,
\end{align}
with the equation of state as 
\begin{equation}
\label{eq:EquationOfState}
	\Ec = \frac{p}{\gamma - 1} + \frac{\rho\left\| \uv \right\|^2}{2} + \frac{\left\| \Bv\right\|^2}{2}.
\end{equation}
where the total mass $\rho$, the momentum $\rho \uv = (\rho u^1, \rho u^2, \rho u^3)^T$, the energy densities $\Ec$ of the plasma system, and the magnetic field vector  $\Bv = (B^1, B^2, B^3)^T$ are all conserved variables.  The velocity $\uv$ and the total pressure $\ptot = p + \frac{1}{2}\lVert \Bv \rVert^2$, together with the hydrodynamic pressure, which is given by the ideal gas law as 
	\begin{equation}
	p = (\gamma -1)(\Ec - \frac{1}{2}\lVert \Bv \rVert^2 - \frac{1}{2}\rho \lVert \uv \rVert^2),
	\end{equation} 
are derived quantities. $\gamma = 5/3$ is the ideal gas constant and the notation  $\| \cdot \|$ is used for the purpose of the Euclidean vector norm.
The MHD equations \eqref{eq:MHD}-\eqref{eq:DivergenceFree} are derived and discussed in many standard plasma textbooks  (e.g., \cite{parks1991physics}).

\subsection{Hyperbolicity of the governing equations}  
The system \eqref{eq:MHD} along with the equation of state \eqref{eq:EquationOfState} comprise a system of hyperbolic conservation laws 
\begin{equation}
\label{eq:hyperbolic_system}
	\qv_{t} +  \nabla \cdot \Fv(\qv) = 0.
\end{equation} 
This set of equations describes the time evolution of all eight conserved variables, $\qv = (\rho, \rho \uv, \Ec, \Bv)$. 
We will denote the Jacobian of the hyperbolic system as $M(\qv) =\partial \Fv/\partial \qv$, which is a diagonalizable matrix with real eigenvalues.
The eigenvalues of the Jacobian matrix represent the wave speeds of the eight waves in MHD system. In particular, in some arbitrary direction  $\nv$  $( \lVert \nv \rVert = 1)$, they can be written as 
\begin{subequations}
\label{eq:eigenvalues}
	\begin{align}
	 & \lambda^{1,8}  = \uv \cdot \nv \mp c_{f}    &\text{fast magnetosonic waves}, \\
	 & \lambda^{2,7}  = \uv \cdot \nv \mp c_{a}     & \text{Alv$\acute{e}$n waves},\\
	 & \lambda^{3,6}  = \uv \cdot \nv \mp c_{s}    & \text{slow magnetosonic waves},\\
	 & \lambda^{4}  = \uv \cdot \nv     & \text{entropy waves}, \\
	 & \lambda^{5}  = \uv \cdot \nv     & \text{divergence waves},
	\end{align} 
\end{subequations}
where
\begin{subequations}
\label{eq:constants}
	\begin{align}
	 & a  \equiv \sqrt{\frac{\gamma p}{\rho}}\\
	 & c_{a}  \equiv \sqrt{\frac{(\Bv \cdot \nv)^2}{\rho}}\\
	 & c_{f}  \equiv       
	\left \{
        \frac{1}{2} \Bigg[a^2 + \frac{\lVert \Bv \rVert^2}{\rho} + \sqrt{\bigg(a^2 + \frac{\lVert \Bv \rVert^2}{\rho} \bigg)^2 - 4a^2\frac{(\Bv \cdot \nv)^2}{\rho}}\Bigg]
        \right \}^{1/2}\\
	 & c_{s}  \equiv  \Bigg\{ \frac{1}{2} \Bigg[a^2 + \frac{\lVert \Bv \rVert^2}{\rho} - \sqrt{\bigg(a^2 + \frac{\lVert \Bv \rVert^2}{\rho} \bigg)^2 - 4a^2\frac{(\Bv \cdot \nv)^2}{\rho}}\Bigg] \Bigg\}^{1/2}
	\end{align}
\end{subequations} 
The eight eigenvalues are well ordered as 
	\begin{equation}
	\lambda ^1 \leq \lambda ^2 \leq \lambda ^3 \leq \lambda ^4 \leq \lambda ^5 \leq \lambda ^6 \leq \lambda ^7 \leq \lambda ^8
	\end{equation}
where the fast and slow magnetosonic waves are nonlinear while the rest of the waves are linearly degenerate.

\subsection{Magnetic potential in 3D}  
\label{sec:3DMP}

There have been many numerical methods presented in the literature for numerically solving the MHD system, but they have faced the main challenge of satisfying the divergence free condition on the magnetic field. Here, we will derive magnetic vector potential equation from the magnetic field equation given in \eqref{eq:MHD}-\eqref{eq:DivergenceFree} system. This will serve as the foundation of our constrained transport framework. 

Since the magnetic field is divergence free, it can always be written as the curl of a magnetic vector potential 
\begin{equation}
	\Bv = \nabla \times \Av.
\end{equation}
\noindent
The key step of the constrained transport scheme is to solve the magnetic potential for correcting the magnetic field. The evolution of the magnetic field in \eqref{eq:MHD} can be written in the following form 
\begin{equation}
	\label{eq:equation20}
	\Bv_{t} + \nabla \times (\Bv \times \uv) = 0,
\end{equation}
using the relation 
\begin{equation*}
	\nabla \cdot (\uv \otimes \Bv - \Bv \otimes \uv) = \nabla \times (\Bv \times \uv). 
\end{equation*}
Since $\Bv$ is divergence free, we set $\Bv = \nabla \times \Av$ and rewrite the evolution equation \eqref{eq:equation20} as 
\begin{equation*}
	\label{eq:equation21}
	\nabla \times \big \{ \Av_{t} + (\nabla\times \Av)\times \uv \big \} = 0.
\end{equation*}
This implies that there exists a scalar function $\psi$ such that
\begin{equation*}
	\Av_{t} + (\nabla \times \Av) \times \uv = - \nabla \psi.
\end{equation*}
There are various choices of the gauge conditions depending on how we chose the $\psi$. Helzel et al. \cite{helzel2011unstaggered} showed that using the Weyl gauge, i.e., setting  $\psi \equiv 0$, one can achieve stable solutions. This condition results in the evolution equation for the magnetic vector potential as 
\begin{equation}
\label{eq:AEquation}
\partial_t \Av + ( \nabla \times \Av ) \times \uv = 0,
\end{equation}
which can be written as:
\begin{equation}
\label{eq:3dmagneticP}
\partial_t \Av + N_1 \Av_x + N_2 \Av_y + N_3 \Av_z = 0,
\end{equation}
with
\begin{align*}
N_1=\begin{bmatrix}
0 & -u^2 & -u^3\\
0 & u^1 & 0 \\
0 & 0 & u^1 \\
\end{bmatrix}, \quad 
N_2=\begin{bmatrix}
u^2 & 0 & 0 \\
-u^1 & 0 & -u^3\\
0 & 0 & u^2\\
\end{bmatrix}, \quad
N_3=\begin{bmatrix}
u^3 & 0 & 0 \\
0 & u^3 & 0 \\
-u^1 & -u^2 & 0\\
\end{bmatrix}.
\end{align*}
The resulting system is only weakly hyperbolic since the matrix of right eigenvectors of the flux Jacobian doesn't have full rank in certain directions. We begin with the flux Jacobian matrix in some arbitrary direction $\nv = (n^1, n^2, n^3) $ to show the weakly hyperbolicity:
\begin{equation}
M(\nv,\uv) : = n^1N_{1} + n^2N_{2} + n^3N_{3} = \begin{bmatrix}
n^2u^2+n^3u^3 & -n^1u^2 & -n^1u^3\\
-n^2u^1 & n^1u^1+n^3u^3 & - n^2u^3 \\
-n^3u^1 & -n^3u^2 & n^1u^1+n^2u^2 \\
\end{bmatrix},
\end{equation}
which has real eigenvalues for all $\lVert \nv \rVert = 1$ as 
\begin{equation}
\lambda ^1 = 0, \quad\quad \lambda^2 = \lambda^3 = \nv \cdot \uv,
\end{equation}
and the right eigenvectors matrix is 
\begin{equation}
R=\Bigg[ r^{(1)} \Big| r^{(2)} \Big| r^{(3)} \Bigg]= \begin{bmatrix}
n^1 & n^2u^3-n^3u^2 & u^1(\uv\cdot \nv) - n^1\lVert \nv \rVert ^2 \\
n^2 & n^3u^1-n^1u^3& u^2(\uv\cdot \nv) - n^2\lVert \nv \rVert ^2  \\
n^3 & n^1u^2-n^2u^1 & u^3(\uv\cdot \nv) - n^3\lVert \nv \rVert ^2  \\
\end{bmatrix}.
\end{equation}
If we assume that $\lVert \uv \rVert \neq 0 $ and $\lVert \nv \rVert = 1$, then the determinant of the matrix of the right eigenvector $R$  can be found as 
\begin{equation*}
	det(R) = - \lVert \uv \rVert ^3 \cos (\alpha) \sin^2(\alpha)
\end{equation*}
where $\alpha$ is the angle between the vectors $\nv$ and $\uv$. Hence, given any nonzero velocity vector $\uv$, there exist four degenerate directions $\alpha = 0, \pi/2, \pi, 3\pi/2$, where the eigenvectors are incomplete due to $det(R) = 0$. Therefore, the system \eqref{eq:3dmagneticP} is only weak hyperbolic.

\subsection{Magnetic potential in 2D}  

In two dimensional case (e.g., in the xy-plane), each conserved variable, $\qv = (\rho, \rho \uv, \Ec, \Bv)$, depends only on $t$, $x$, and $y$ independent variables. Thus, the divergence free condition can be simplified as 
\begin{equation*}
	\nabla \cdot \Bv = B^{1}_{x} + B^{2}_{y}.
\end{equation*}
There is no point to define $B^{3}$ since any value of $B^3$ satisfies the divergence free condition in 2D. The evolution of the magnetic field involves only the third component of the magnetic vector potential, so that magnetic field components, $B^1$ and $B^2$, can be found as 
\begin{equation}
	\label{eq:2drelation}
	B^1 = A^{3}_{y} \quad \text{and}  \quad B^2 = - A^{3}_{x}.
\end{equation}
Since we only need the third component of the magnetic vector potential to define magnetic field, the constrained transport scheme in 2D is reduced to the magnetic vector potential to a magnetic scalar potential equation as 
 \begin{equation}
	\label{eq:2dmagneticP}
 	A^{3}_{t}+ u^{1} A^{3}_{x} + u^{2} A^{3}_{y} = 0
\end{equation} 
This scalar advection equation is strongly hyperbolic by contrast with the 3D case. 

\section{Outline of the constrained transport methodology}
\label{sec:CT}
 The major challenge when dealing with the numerical solution of the system \eqref{eq:MHD}-\eqref{eq:EquationOfState} is to satisfy the divergence-free condition \eqref{eq:DivergenceFree}. The constrained transport methodology is one of the dominant approach in the literature to overcome this challenge. In this section we will give the main idea of this method and our perspective how we are using high order CT framework.
The idea of the constrained transport methodology is based on updating the conserved variables $\qv = (\rho, \rho \uv, \Ec, \Bv)$. One of the early works on this idea of using the magnetic vector potential equations for CT solution to the MHD equations is presented by Wilson in \cite{wilson1975some}, as well as Dorfi \cite{dorfi1986numerical}. However, modern shock capturing strategies were not used in those works and hence caused strong numerical diffusion. Londrillo and Del Zanna \cite{londrillo2000high} used the magnetic potential solutions in the context of shock capturing methods, along with De Sterck \cite{de2001multi}, Londrillo and Del Zanna \cite{londrillo2004divergence}, and Rossmanith \cite{rossmanith2006unstaggered}. Helzet et al. \cite{helzel2013high} developed a framework for unstaggered CT scheme coupling a conservative finite volume hyperbolic scheme for the MHD equations with a non conservative finite volume scheme for the magnetic vector potential equation. In our previous work \cite{christlieb2014finite}, we developed a high order finite difference constrained transport scheme. In this paper, we present a kernel based high order unconditionally stable CT method. We correct the conserved variables, $\qv$, using the magnetic vector potential formulation with respect to relation $\Bv = \nabla \times \Av$. Here, we give an outline of the general unstaggered CT framework listing all important steps. 
Consider a semi discrete system of ordinary differential equations for MHD equations \eqref{eq:MHD} 
\begin{equation}
	\qv'_{mhd}(t) = \mathcal{L}_{1}(\qv_{mhd}(t))
\end{equation}
where $\qv_{mhd}(t)$ represents the grid function at time $t$ consisting of all point-wise values
 of the conserved quantities in the ideal MHD system $\qv_{mhd} = (\rho, \rho \uv, \Ec, \Bv)$. The details about $\mathcal{L}_{1}(\qv_{mhd}(t))$ were presented in \cite{christlieb2014finite}.

And also consider the update for the magnetic potential equation (\eqref{eq:2dmagneticP} for 2D case and \eqref{eq:3dmagneticP} for 3D case) on the mesh,
 \begin{equation}
	\qv'_{A}(t) = \mathcal{L}_{2}(\qv_{A}(t), \qv_{mhd}(t)),
\end{equation}
where $\qv_{A}$ denotes vector potential $\Av$ in 3D or the scalar potential $A^3$ in 2D.

The key steps advancing the solution from its current time step $t = t^n$ (or the initial condition at $t^0$) to its new time step $t^{n+1}$ are listed below: 

\begin{itemize}
\item step 0: Start with the given current time step $\qv_{mhd}^{n} = (\rho^n, \rho \uv^n, \Ec^n, \Bv^n)^T$ and $\qv_{A}^{n}$. 

\item step 1: Obtain $\qv_{mhd}^{*}$ and $\qv_{A}^{n+1}$ separately, where $$\qv_{mhd}^{*}=\left(\rho^{n+1}, \rho^{n+1} \uv^{n+1}, \Ec^{*}, \Bv^{*} \right).$$
Here, $\Ec^{*}$ and $\Bv^{*}$ are given with a $*$ superscript instead of $n+1$ to indicate that the predicted $\Bv$ and $\Ec$ will be corrected by a predictor-corrector Constrained Transpose method before the end of the time step.

\item step 2: Replace $\Bv^*$ to $\Bv^{n+1}$ by a discrete curl of $\qv_{A}^{n+1}$.
$$ \Bv^{n+1} = \nabla \times \qv_{A}^{n+1}. $$

\item step 3: Set the corrected total energy density value $\Ec^{n+1}$ based on one of the following options: 
\par\qquad Option 1: Keep the total energy conserved 
$$ \Ec^{n+1} = \Ec^*.$$
\par\qquad Option 2: Keep the pressure the same after updating the magnetic field $$ \Ec^{n+1} = \Ec^* + \frac{1}{2}(\left\| \Bv^{n+1}\right\|^2 - \left\| \Bv^*\right\|^2)。$$
(Second option sometimes helps to prevent negative pressure).
\end{itemize}
In Section 5, we show that our approach leads to a divergence free solution.
We now describe how we constraint the update for $\qv_{A}^{n+1}$ using our kernel based approach.    
\section{Hamilton-Jacobi equations}

In this section we introduce the main ideas in our kernel-based method \cite{christlieb2017kernel, christlieb2019kernel}, which is derived from the Method of Lines Transpose (MOL$^T$) \cite{causley2013method, causley2014method, causley2014higher, causley2016method, christlieb2016weno, causley2017method}. The simplest way to describe MOLT is: we start by discretizing the problem in time; we then use  a global approximation for the inherently local term (the derivative).  By doing this, we are able to make an explicit approximation unconditionally stable.  As for the global approximation, think FFT, but with a different convolution kernel.  The form of the approximation is what facilitates the stability of the method.  The method is an O(N) because the convolution with the kernels can be evaluated using a three term recreation. 

\subsection{1D Hamilton-Jacobi equations} 
\label{sec:1dformulation}

Consider a 1D Hamilton-Jacobi equation
\begin{equation}
\label{eq:HJ}
A_{t}+H(A_{x})=0,\qquad A(x,0) = A^{0}(x), \qquad x\in[a,b],
\end{equation}
where $H$ is the Hamiltonian flux. 
For the time discretization purpose to evolve the solution from time $t^{n}$ to $t^{n+1}$, we use the classical explicit SSP RK schemes \cite{gottlieb2001strong}. In this work, we propose to use the following SSP RK schemes such as the first order forward Euler scheme 
\begin{align}
\label{eq:rk1}
A^{n+1}=A^{n}-\Delta t \hat{H}(A^{n,-}_{x},A^{n,+}_{x});
\end{align} 
the second order SSP RK scheme
\begin{align}
\label{eq:rk2}
& A^{(1)}=A^{n}-\Delta t \hat{H}(A^{n,-}_{x},A^{n,+}_{x}),\nonumber\\
& A^{n+1}=\frac{1}{2}A^{n}+\frac{1}{2}\left( A^{(1)} -\Delta t \hat{H}(A^{(1),-}_{x},A^{(1),+}_{x}) \right);
\end{align}
and the third order SSP RK scheme
\begin{align}
\label{eq:rk3}
& A^{(1)}=u^{n}-\Delta t \hat{H}(A^{n,-}_{x},A^{n,+}_{x}),\nonumber\\
& A^{(2)}=\frac{3}{4}A^{n}+\frac{1}{4} \left( A^{(1)}-\Delta t \hat{H}(A^{(1),-}_{x},A^{(1),+}_{x}) \right), \nonumber\\
& A^{n+1}=\frac{1}{3}A^{n}+\frac{2}{3} \left( A^{(2)}-\Delta t \hat{H}(A^{(2),-}_{x},A^{(2),+}_{x}) \right).
\end{align}
where $\Delta t$ denotes the time step and $\hat{H}$ is the the numerical Hamiltonian, e.g., the local Lax-Friedrichs Hamiltonian flux \cite{jiang2000weighted}:
		\begin{equation}
		\hat{H}(\phi^{-},\phi^{+})=H(\frac{\phi^{-}+\phi^{+}}{2}) -c(\phi^{-},\phi^{+})\frac{(\phi^{+}-\phi^{-})}{2},
		\end{equation}
	with $c(\phi^{-},\phi^{+})=\max_{\phi\in[\min(\phi^{-},\phi^{+}), \max(\phi^{-},\phi^{+})]} |H'(\phi)|$.
Here, $A^{-}_{x}$ and $A^{+}_{x}$ are one-side derivatives with \emph{left}-biased and \emph{right}-biased methods, respectively, to approximate $A_{x}$. Finding these derivatives has the dominant role in our work and we will present the details of the construction of them in the following subsection. 


\subsection{Approximation of the first order derivative $\partial_x$}  

Approximating the partial spatial derivative terms with kernel based scheme is the major part of this work. In this section we briefly review the construction of the $\partial_x$ derivative approximation using kernel based formulation established in \cite{christlieb2017kernel}.
Let's define operators $\mathcal{L}_{L}$ and $\mathcal{L}_{R}$, and their inverse operators
\begin{subequations}
\label{eq:DL}
	\begin{align}
		& \mathcal{L}_{L}=\mathcal{I}+\frac{1}{\alpha}\partial_{x} \ \Rightarrow \
		\mathcal{L}_{L}^{-1}[v,\alpha](x)
		=\alpha\int_{a}^{x}e^{-\alpha(x-y)}v(y)dy +A_{L}e^{-\alpha(x-a)}, \\
		& \mathcal{L}_{R}=\mathcal{I}-\frac{1}{\alpha}\partial_{x} \ \Rightarrow \
		\mathcal{L}_{R}^{-1}[v,\alpha](x)
		=\alpha\int_{x}^{b}e^{-\alpha(y-x)}v(y)dy +B_{R}e^{-\alpha(b-x)}, 
	\end{align}
\end{subequations}
where $\mathcal{I}$ is the identity operator and $\alpha$ is a positive constant. We use $I^{L}$ and  $I^{R}$ to denote the convolution integral as 
\begin{equation}
	I^{L} = \alpha\int_{a}^{x}e^{-\alpha(x-y)}v(y)dy, \quad
	I^{R} = \alpha\int_{x}^{b}e^{-\alpha(y-x)}v(y)dy. 
\end{equation}
We can see that $I^{L}$ depends on the function values of $v$ from left end point $a$ to $x$, as well as $I^{R}$ is for right end point $b$ to $x$. Also note that $A_{L}$ and $B_{R}$ are determined by the boundary conditions. For instance, if we assume periodic boundary conditions, i.e.,
\begin{equation}
	\mathcal{L}_{L}^{-1}[v,\alpha](a) = \mathcal{L}_{L}^{-1}[v,\alpha](b), \quad \text{and}, \quad  \mathcal{L}_{R}^{-1}[v,\alpha](a) = \mathcal{L}_{R}^{-1}[v,\alpha](b).
\end{equation}
then we obtain 
\begin{equation}
	A_{L} = \frac{I^{L}[v,\alpha](b)}{1-\mu}, \quad \text{and}  \quad  B_{R} = \frac{I^{R}[v,\alpha](a)}{1-\mu},
\end{equation}
where $\mu = e^{-\alpha(b-a)}$. 

Now, let's define new operators $\mathcal{D}_{L}$ and $\mathcal{D}_{R}$, 
		\begin{align*}
		\mathcal{D}_{L}=\mathcal{I}-\mathcal{L}_{L}^{-1}, 
		\quad \text{and} \quad 
		\mathcal{D}_{R}=\mathcal{I}-\mathcal{L}_{R}^{-1}.
		\end{align*} 
Then $\frac{1}{\alpha}\partial_{x}$ can be represented using the infinite series of these operators as
\begin{subequations}
	\begin{align}
	& \frac{1}{\alpha}\partial_{x} = \mathcal{L}_{L} - \mathcal{I} =  \mathcal{L}_{L}( \mathcal{I} - \mathcal{L}_{L}^{-1})  ={\mathcal{D}_{L}}{(\mathcal{I} - \mathcal{D}_{L})^{-1}}=  \sum_{p=1}^{\infty}\mathcal{D}_{L}^{p}\\
	& \frac{1}{\alpha}\partial_{x} = \mathcal{I} - \mathcal{L}_{R}  =  \mathcal{L}_{R}( \mathcal{L}_{R}^{-1} - \mathcal{I}) =  {-\mathcal{D}_{R}}{(\mathcal{I} - \mathcal{D}_{R})^{-1}} = - \sum_{p=1}^{\infty}\mathcal{D}_{R}^{p}
	\end{align}
\end{subequations}	
Hence, if $A$ is a periodic function, then we can obtain the approximations of the first derivatives $A_{x}^{\pm}$ using partial sums as
\begin{subequations} 
\label{eq:partialsum_per}  
	\begin{align}
	A_{x}^{-}
	\approx \left\{ \begin{array}{ll}
	 \alpha\sum_{p=1}^{k}\mathcal{D}_{L}^{p}[A,\alpha](x), & k=1,2, \\
	 \\
	 \alpha\sum_{p=1}^{k}\mathcal{D}_{L}^{p}[A,\alpha](x) - \alpha \mathcal{D}_{0}*\mathcal{D}^2_{L}[A,\alpha](x), & k=3,\\
	 \end{array}
	 \right.
	\end{align}  
\text{and}
	\begin{align}    
	A_{x}^{+}
	\approx \left\{ \begin{array}{ll}
	 -\alpha\sum_{p=1}^{k}\mathcal{D}_{R}^{p}[A,\alpha](x), & k=1,2, \\
	 \\
	- \alpha\sum_{p=1}^{k}\mathcal{D}_{R}^{p}[A,\alpha](x) + \alpha \mathcal{D}_{0}*\mathcal{D}^2_{R}[A,\alpha](x), & k=3.\\
	 \end{array}
	 \right.
	\end{align}
\end{subequations} 
For $k=3$ case, there is an additional term $\mathcal{D}_{0}$ which is defined as 
	\begin{equation}
	\mathcal{D}_{0}[v,a] = v(a) - \frac{\alpha}{2}\int_{a}^{b}e^{-\alpha|x-y|}v(y)dy - A_{0}e^{-\alpha(x-a)} - B_{0}e^{-\alpha(b-x)}.
	\end{equation}
where $A_{0}$ and $B_{0}$ are determined by boundary condition. For example, if $\mathcal{D}_{0}$ is a periodic function such that
$$\mD_{0}[v,\alpha](a)=\mD_{0}[v,\alpha](b), \quad\text{and}\quad \partial_{x}\mD_{0}[v,\alpha](a)=\partial_{x}\mD_{0}[v,\alpha](b),$$
then we get
\begin{align}
	\label{eq:gxx_bc_per}
A_{0} = \frac{I^{0}[v,\alpha](b)}{1-\mu}, \quad
B_{0} = \frac{I_{0}[v,\alpha](a)}{1-\mu},
\end{align}
with $I^{0}[v,\alpha](x)=\frac{\alpha}{2}\int_{a}^{b}e^{-\alpha|x-y|}v(y)dy$. 

Here, we take $\alpha = \beta/(c\Delta t)$, where c is the maximum wave speed and $\beta$ is a constant independent of the time step $\Delta t$. Therefore, the accuracy for the approximation \eqref{eq:partialsum_per} to the $\partial_{x}A$ is $\mathcal{O}(\Delta t^k)$.
If $\beta$ is chosen appropriately, then the scheme is A-stable and hence allows for large time step evolution. Additionally, the linear stability of the method which argued in the following theorem has been proven in \cite{christlieb2017kernel}. 
\begin{thm}\label{thm4}
For the linear equation $\phi_{t}+c\phi_{x}=0$, (i.e. the Hamiltonian is linear) with periodic boundary conditions, we consider the $k^{th}$ order SSP RK method as well as the  $k^{th}$ partial sum in \eqref{eq:partialsum_per}, with $\alpha=\beta/(|c|\Delta t)$. Then there exists a constant $\beta_{k,max}>0$ for $k=1,\,2,\,3$, such that the scheme is A-stable provided  $0<\beta\leq\beta_{k,\max}$. The constants $\beta_{k,max}$ for $k=1,\,2,\,3$ are summarized in Table \ref{tab0}.
\end{thm}
\begin{table}[htb]
	\caption{\label{tab0}\em $\beta_{k,\max}$ in Theorem \ref{thm4} for  $k=1,\,2,\,3$.}
	\centering
	\vspace{0.3cm}
	\begin{tabular}{| l | p{1cm} | p{1cm} | p{1cm} |}
		\hline
		$k$ &  1  & 2  & 3  \\\hline
		$\beta_{k,max}$  &  2  &  1  &  1.243  \\\hline
	\end{tabular}
\end{table}
	


\subsection{Non periodic boundary conditions }  

Based on the boundary conditions, we need to modify the treatments for $D^p_{L}$, $D^p_{R}$, and $D_{0}$ for non periodic boundary conditions so that the boundary conditions specified on these operators are consistent with the boundary conditions imposed on $A$. 


Firstly, we examine the operators $D_{*}$ for non periodic boundary condition, where $*$ denotes $0$, $L$, or $R$. Assume that we are given numbers $C_{a}$ and $C_{b}$.
\begin{itemize}
	\item If we require
	$$\mD_{L}[v,\alpha](a)=C_{a}, \quad\text{and}\quad \mD_{R}[v,\alpha](b)=C_{b},$$
	then, the boundary coefficients are obtained as
	\begin{align}
	\label{eq:fx_bc_dir}
	A_{L}=v(a) - C_{a} , \quad \text{and} \quad B_{R}= v(b) - C_{b}.
	\end{align}
	\item If we require
	$$\mD_{0}[v,\alpha](a)=C_{a}, \quad\text{and}\quad \mD_{0}[v,\alpha](b)=C_{b},$$
	then, the boundary coefficients are obtained as
	\begin{subequations}
		\label{eq:gxx_bc_dir}
		\begin{align}
		& A_{0}=\frac{1}{1-\mu^2}\left( \mu\left(I^{0}[v,\alpha](b)-v(b)+C_{b}\right) - \left(I^{0}[v,\alpha](a)-v(a)+C_{a}\right)\right), \\
		& B_{0}=\frac{1}{1-\mu^2}\left( \mu\left(I^{0}[v,\alpha](a)-v(a)+C_{a}\right) - \left(I^{0}[v,\alpha](b)-v(b)+C_{b}\right)\right).
		\end{align}
	\end{subequations}
\end{itemize}

%

Next, we introduce a modification of the partial sums \eqref{eq:partialsum_per} to achieve a higher order accuracy for the non-periodic case.
Assume that we have some derivative values at boundary, i.e., $\partial_{x}^{m}A(a)$ and $\partial_{x}^{m}A(b)$, $m\geq1$.
Using integration by parts one can derive the following modified partial sums for $k\leq3$ to deal with the non-periodic boundary conditions 
\begin{subequations}
	\label{eq:partialsum_dir}
	\begin{align}
	A_{x}^{-}(x)\approx \mathcal{P}^{-}_{k}[A,\alpha] = \left\{\begin{array}{ll}
	\alpha\sum\limits_{p=1}^{k}\mathcal{D}_{L}[A_{1,p},\alpha](x), & k=1,\, 2,\\
	\alpha\sum\limits_{p=1}^{k}\mathcal{D}_{L}[A_{1,p},\alpha](x) -\alpha \mD_{0}[A_{1,3},\alpha](x), & k=3,\\
	\end{array}
	\right.
	\end{align}
	\begin{align}
	A_{x}^{+}(x)\approx \mathcal{P}^{+}_{k}[A,\alpha] = \left\{\begin{array}{ll} -\alpha\sum\limits_{p=1}^{k}\mathcal{D}_{R}[A_{2,p},\alpha](x),& k=1,\, 2,\\
	-\alpha\sum\limits_{p=1}^{k}\mathcal{D}_{R}[A_{2,p},\alpha](x) +\alpha \mD_{0}[A_{2,3},\alpha](x), & k=3.\\
	\end{array}
	\right.
	\end{align}
\end{subequations}
And $A_{1,p}$ and $A_{2,p}$ are given as
\begin{subequations}
	\label{eq:expression}
	\begin{align}
	& \left\{\begin{array}{ll}
	A_{1,1}=A,\\
	\displaystyle A_{1,2}=\mathcal{D}_{L}[A_{1,1},\alpha] - \sum_{m=2}^{k}\left(-\frac{1}{\alpha}\right)^{m} \partial_{x}^{m}A(a) e^{-\alpha(x-a)},\\
	\displaystyle A_{1,3}=\mathcal{D}_{L}[A_{1,2},\alpha] + \sum_{m=2}^{k}(m-1)\left(-\frac{1}{\alpha}\right)^{m} \partial_{x}^{m}A(a) e^{-\alpha(x-a)},\\
	\end{array}
	\right.\\
	& \left\{\begin{array}{ll}
	A_{2,1}=A, \\
	\displaystyle A_{2,2}=\mathcal{D}_{R}[A_{2,1},\alpha] - \sum_{m=2}^{k}\left(\frac{1}{\alpha}\right)^{m} \partial_{x}^{m}A(b) e^{-\alpha(b-x)}, \\
	\displaystyle A_{2,3}=\mathcal{D}_{R}[A_{2,2},\alpha] + \sum_{m=2}^{k}(m-1)\left(\frac{1}{\alpha}\right)^{m} \partial_{x}^{m}A(b) e^{-\alpha(b-x)},\\
	\end{array}
	\right.
	\end{align}
\end{subequations}
where the boundary conditions for the operators are imposed as
	\begin{align*}
	& \alpha\mathcal{D}_{L}[A_{1,1},\alpha](a)=A_{x}(a), \quad
	\alpha\mathcal{D}_{R}[A_{2,1},\alpha](b)=-A_{x}(b),\\
	& \alpha\mathcal{D}_{L}[A_{1,p},\alpha](a)=\alpha\mathcal{D}_{R}[A_{2,p},\alpha](b)= 0,
	\quad \text{for} \ p\geq2,\\
	& \alpha\mD_{0}[A_{*,3},\alpha](a) = \alpha\mD_{0}[A_{*,3},\alpha](b)=0, \quad \text{$*$ could be 1 or 2.}
	\end{align*} 
The modified partial sum \eqref{eq:partialsum_dir} agrees with the derivative values at the boundary
$$ \mathcal{P}^{-}_{k}[A,\alpha](a)=A_{x}(a), \quad
\mathcal{P}^{+}_{k}[A,\alpha](b)=A_{x}(b). $$
Furthermore, we have the following theorem, which is a result of the Theorem 2.3 from \cite{christlieb2019kernel}. 
\begin{thm}
	Suppose $A\in\mathcal{C}^{k+1}[a,b]$. If we take $\alpha=\beta/(c\Delta t)$, then the modified partial sums \eqref{eq:partialsum_dir} satisfy
	\begin{align}
			\| A_{x} - \mathcal{P}^{-}_{k}[A,\alpha] \|_{\infty} = \mathcal{O}(\Delta t^k), \quad
			\| A_{x} - \mathcal{P}^{+}_{k}[A,\alpha] \|_{\infty}  = \mathcal{O}(\Delta t^k), \quad k=1, \, 2,\, 3.
	\end{align}
\end{thm}

\par As we have seen in \eqref{eq:expression}, we need the derivatives $\partial_{x}^{m}A(a)$ and $\partial_{x}^{m}A(b)$, $m\geq1$, for non periodic boundary conditions. Here, we will only focus on the outflow boundary conditions. High order extrapolations are used to get all the derivatives at boundaries (that is, there are no physical boundary given). The details of the general non-periodic conditions can be found in our previous work \cite{christlieb2019kernel}.

\subsection{Space discretization}

According to the MOL$^T$ method, a boundary value problem is obtained by discretizing the equation in time first and an integral formulation is presented to solve the boundary value problem.  A fully discrete numerical solution is then obtained by discretizing $\mathcal{D}_{L}$ and $\mathcal{D}_{R}$ operators in space. In this subsection we will give the details of the spatial discretization of the  $\mathcal{D}_{L}$ and $\mathcal{D}_{R}$ operators and WENO-based quadrature formulation to approximate the convolution integrals appear in the  $\mathcal{D}_{L}$ and $\mathcal{D}_{R}$ operators. 

Suppose we divide the domain $[a,b]$ with $N+1$ uniformly distributed grid points
$$\Delta x = (b-a)/N, \quad x_{i}=a+i \Delta x, \quad i=0, 1,\ldots, N.$$
The convolution integrals $I^{L}_{i}=I^L[v,\alpha](x_i)$ and $I^{R}_{i}=I^R[v,\alpha](x_i)$ satisfy a recursive relation 
\begin{subequations}
	\label{eq:recursive}
	\begin{align}
	& I^L_i = e^{-\alpha\Delta x_{i}} I^L_{i-1} + J^L_i,\quad i=1,\ldots,N, \quad I^L_0 = 0, \\
	& I^R_i = e^{-\alpha\Delta x_{i+1}} I^R_{i+1} + J^R_i,\quad i=0,\ldots,N-1, \quad I^R_N = 0,
	\end{align}
\end{subequations}
where,
\begin{align}
\label{eq:JLR}
J^L_{i} =  \alpha \int_{x_{i-1}}^{x_{i}} v(y)e^{-\alpha (x_{i}-y)}dy,\ \ \ \
J^R_{i} =  \alpha \int_{x_{i}}^{x_{i+1}} v(y)e^{-\alpha (y-x_{i})}dy.
\end{align}

In \cite{christlieb2017kernel}, we have presented a high order and robust framework to calculate $J^L_{i}$ and $J^R_{i}$. Here we will provide a brief description of methodology. 

Note that, to approximate $J^L_{i}$ with $k^{th}$ order accuracy, we may choose the interpolation stencil 
$S(i)=\left\{x_{i-r},\ldots,x_{i-r+k} \right\}$,
which contains $x_{i-1}$ and $x_{i}$. There is a unique polynomial $p(x)$ of degree at most $k$ that interpolates $v(x)$ at the nodes in $S(i)$. Then $J^L_{i}$ is approximated by
\begin{align}
\label{eq:linear_JL}
J^L_{i}\approx\alpha \int_{x_{i-1}}^{x_{i}} p(y)e^{-\alpha (x_{i}-y)}dy.
\end{align}
Note that, the integral on the right hand side can be evaluated exactly.
Similarly, we can approximate $J^R_{i}$ by
\begin{align}
\label{eq:linear_JR}
J^R_{i}\approx\alpha \int_{x_{i}}^{x_{i+1}} p(y)e^{-\alpha (y-x_{i})}dy,
\end{align}
with polynomial $p(x)$ interpolating $v(x)$ on stencil $S(i)=\left\{x_{i+r-k},\ldots,x_{i+r} \right\}$, which includes $x_{i}$ and $x_{i+1}$.

Since the quadratures with a fixed stencil for approximating the $J^{L}_{i}$ and $J^{L}_{i}$ may develop spurious oscillations which violates the entropy solutions, we use WENO-based quadrature formula and the nonlinear filter to control oscillations and capture the correct solution. Such a methodology is discussed in \cite{christlieb2017kernel} and we will modify that method to solve HJ equations since both equations involve discontinuous derivatives for the solutions and hence, generate the entropy violating solutions.

To summarize, if $A_{x}$ is periodic function, then we will be using the following modified sums framework from \cite{christlieb2019kernel} instead of \eqref{eq:partialsum_per} for approximation of $A_{x}^{\pm}$ at $x_i$:
\begin{subequations}
\label{eq:change_per}
\begin{align}
& A_{x,i}^{-}=\alpha\mathcal{D}_{L}[A,\alpha](x_{i}) + \alpha\sum_{p=2}^{k}\sigma_{i,L}^{p-2}\mathcal{D}_{L}^{p}[A,\alpha](x_{i}),\\
& A_{x,i}^{+}=-\alpha\mathcal{D}_{R}[A,\alpha](x_{i}) - \alpha\sum_{p=2}^{k}\sigma_{i,R}^{p-2}\mathcal{D}_{R}^{p}[A,\alpha](x_{i});
\end{align}
\end{subequations}
and if $A_{x}$ is non periodic function, then we we will use the following formulation instead of \eqref{eq:partialsum_dir}
\begin{subequations}
	\label{eq:change_dir}
	\begin{align}
	& A_{x,i}^{-}=\alpha\mathcal{D}_{L}[A_{1,1},\alpha](x_{i}) + \alpha\sum_{p=2}^{k}\sigma_{i,L}^{p-2}\mathcal{D}_{L}[A_{1,p},\alpha](x_{i}),\\
	& A_{x,i}^{+}=-\alpha\mathcal{D}_{R}[A_{2,1},\alpha](x_{i}) - \alpha\sum_{p=2}^{k}\sigma_{i,R}^{p-2}\mathcal{D}_{R}[A_{2,p},\alpha](x_{i}).
	\end{align}
\end{subequations} 
We only use the WENO formulation when  $p = 1$ while we use cheap high order linear formulation for the case $p > 1$. The filters $\sigma_{i,L}$ and $\sigma_{i,R}$ are obtained based on the smoothness indicators from the WENO quadrature.

\begin{figure}
	\centering
	\includegraphics[width=0.5\textwidth]{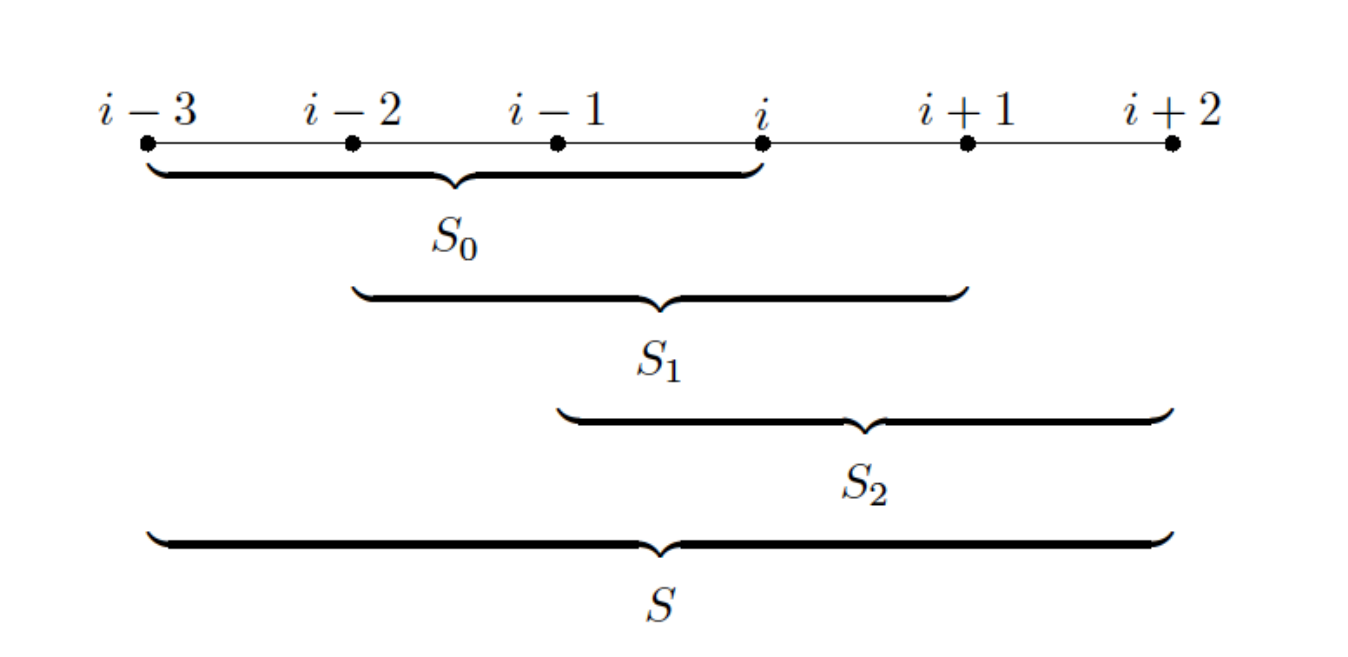}
	\caption{\em The structure of the stencils in WENO integration.   }
	\label{Fig0}
\end{figure}

The WENO quadrature for  $J^{L}_{i}$ is presented here as an example. The related stencil is given in Figure \ref{Fig0}. We choose the big stencil as $S(i)=\{x_{i-3},\ldots, x_{i+2}\}$ and the three small stencils as $S_{r}(i) =\{ x_{i-3+r},\ldots, x_{i+r}\}$, $r = 0, 1, 2$. We will only present the formulas for the case of a uniform mesh, i.e., $\Delta x_{i}=\Delta x$ for all $i$ here, but the WENO methodology is still applicable to the case of a nonuniform mesh, see \cite{shu2009high}. 

\begin{enumerate}
	\item
	We approximate the integrals on each small stencils $S_{r}(i)$ as follows
	\begin{equation}
	\label{eq:weno1}
	J^{L}_{i,r} =\alpha \int_{x_{i-1}}^{x_{i}}e^{-\alpha(x_{i}-y)}p_{r}(y)dx,
	\end{equation}
	where $p_{r}(x)$ is the polynomial interpolating $v(x)$ on nodes $S_{r}(i)$.
	
	\item
	Similarly, on the big stencil $S(i)$, we obtain
	\begin{equation}
	\label{eq:weno2}
	J^{L}_{i} =\alpha \int_{x_{i-1}}^{x_{i}}e^{-\alpha (x_{i}-y)}p(y)dx=\sum_{r=0}^{2}d_{r}J^{L}_{i,r},
	\end{equation}
	with the linear weights $d_{r}$ satisfying $\sum_{r=0}^{2}d_{r}=1$.
	
	\item
	We develop the following nonlinear weights $\omega_{r}$ using  the linear weights $d_{r}$ 
	\begin{equation}
	\omega_{r}=\tilde{\omega}_{r}/\sum\limits_{s=0}^{2}\tilde{\omega}_{s}, \ \ r=0,\ 1,\ 2,
	\end{equation}
	with
	\begin{equation*}
	\tilde{\omega}_{r}=d_{r}\left( 1+\frac{\tau_{5}}{\epsilon+\beta_{r}}\right) .
	\end{equation*}
	We  take $\epsilon=10^{-6}$ as a small positive number, $\epsilon>0$,  in our numerical test problems to avoid zero at the denominator. The smoothness indicator $\beta_{r}$ is determined as 
	\begin{equation}
	\beta_{r}=\sum_{l=2}^{3} \int_{x_{i-1}}^{x_{i}}\Delta x_{i}^{2l-3} \left(\frac{\partial^{l}p_{r}(x)}{\partial x^{l}}\right)^2 dx,
	\end{equation}
	which is used to  measure the relative smoothness of the function $v(x)$ in the stencil $S_{r}(i)$. In particular, we have the expressions as 
	\begin{subequations}
		\begin{align*}
		& \beta_{0} = \frac{13}{12} ( -v_{i-3} + 3v_{i-2} - 3v_{i-1} + v_{i} )^2
		+ \frac{1}{4} ( v_{i-3} - 5v_{i-2} + 7v_{i-1} - 3v_{i} )^2,\\
		& \beta_{1} = \frac{13}{12} (-v_{i-2} + 3v_{i-1} - 3v_{i} + v_{i+1} )^2
		+ \frac{1}{4} ( v_{i-2} - v_{i-1} - v_{i} + v_{i+1} )^2,\\
		& \beta_{2} = \frac{13}{12} (-v_{i-1} + 3v_{i} - 3v_{i+1} + v_{i+2} )^2 
		+ \frac{1}{4} (-3v_{i-1} + 7v_{i} - 5v_{i+1} + v_{i+2} )^2.
		\end{align*}
	\end{subequations}
	Here, $\tau_{5}$ states the absolute difference between $\beta_{0}$ and $\beta_{2}$ 
	$$\tau_{5}=|\beta_{0}-\beta_{2}|.$$
	Furthermore, we introduce a parameter $\xi_{i}$ as
	\begin{align}
	\xi_{i}=\frac{\beta_{min}}{\beta_{max}},
	\end{align}
	which will be use to create the nonlinear filter. Here,
	\begin{subequations}
		\begin{align*}
		\beta_{max} = 1 + \left(\frac{\tau_{5} }{  \epsilon + \min(\beta_{0},\beta_{2}) } \right)^2, \quad
		\beta_{min} = 1 + \left(\frac{\tau_{5} }{  \epsilon + \max(\beta_{0},\beta_{2}) } \right)^2.
		\end{align*}
	\end{subequations}
	Note that we have developed the nonlinear weights using the idea of the WENO-Z method proposed in \cite{borges2008improved}, which has less dissipation and higher resolution comparing to the original WENO method.
	\item
	Lastly, we obtain the approximation
	\begin{equation}
	J^{L}_{i}=\sum_{r=0}^{2}\omega_{r}J^{L}_{i,r}.
	\end{equation}
	The filter $\sigma_{i,L}$ is determined as 
	\begin{align}
	\sigma_{i,L}=\min (\xi_{i-1}, \xi_{i}).
	\end{align}
	
\end{enumerate}
The process to obtain $J^{R}_{i}$ and $\sigma_{i,R}$ is mirror symmetric to that of $J^{L}_{i}$ and $\sigma_{i,L}$ with
respect to point $x_{i}$.

\subsection{2D magnetic potential equation}

According to the Constrained Transport formulation described in Section \ref{sec:CT}, we must update the solution of the magnetic potential equation by solving a discrete version of the following equation
	\begin{equation}
	\label{eq:2dmag-pot1}   
	A^{3}_{t} + u^1(x, y)A^{3}_{x} +  u^2(x, y)A^{3}_{y} = 0
	\end{equation}
where the velocity components $u^1$ and $u^2$ are known from the previous time step due to the solution of the base part, HCL. Since the velocity functions are given, we can consider \eqref{eq:2dmag-pot1}  as a Hamilton-Jacobi equation
	\begin{equation}
	\label{eq:2dmag-pot2}   
	A^{3}_{t} + H(A^{3}_{x} , A^{3}_{y}) = 0
	\end{equation}
with Hamiltonian flux
	\begin{equation}
	\label{eq:2dmag-pot3}   
	H(A^{3}_{x} , A^{3}_{y}) = u^1(x, y)A^{3}_{x} +  u^2(x, y)A^{3}_{y}.
	\end{equation}
We can directly apply a two dimensional version of the kernel-based framework presented in Section \ref{sec:1dformulation}. The 2D semi discrete scheme can be written as 
	\begin{equation}
	\label{eq:2dmag-pot4}  
	\frac{dA_{i,j}^{3}(t)}{dt} = - \hat{H}(A^{3-}_{x}|_{i,j} , A^{3+}_{x}|_{i,j}, A^{3-}_{y}|_{i,j} , A^{3+}_{y}|_{i,j}) 
	\end{equation}
 on each point $(x_i, y_j)$, where $\hat{H}$ is a Lipschitz continuous Hamiltonian flux. We use the Lax-Friedrichs flux and the equation \eqref{eq:2dmag-pot4} becomes

	\begin{align}
	\label{eq:2dmag-pot5}  
	\frac{dA_{i,j}^{3}(t)}{dt} =&  -u_{i,j}^1\Big( \frac{A^{3-}_{x}  + A^{3+}_{x}}{2} \Big)|_{i,j} - u_{i,j}^2\Big( \frac{A^{3-}_{y} + A^{3+}_{y}}{2} \Big)|_{i,j} \notag\\
	& + c_1\Big( \frac{A^{3+}_{x} - A^{3-}_{x}}{2} \Big)|_{i,j} + c_2\Big( \frac{A^{3+}_{y} - A^{3-}_{y}}{2} \Big)|_{i,j} 
	\end{align}
with
	\begin{equation*}
	c_1 = \max_{i,j}{|u^1_{i,j}|}  \quad \text{and} \quad c_2 = \max_{i,j}{|u^2_{i,j}|}.
	\end{equation*}
We remark that the scheme \eqref{eq:2dmag-pot5} with this global $c_m$ can be very dissipative for some Hamilton-Jacobi equations. There is another way of choosing $c_m$ which is by taking the max value on the local stencil.
\par The approximations $A^{3\pm}_{x}|_{i,j}$ and $A^{3\pm}_{y}|_{i,j}$ to the derivatives of functions $A_x(x, y)$ and $A_y(x, y)$ at $(x_i, y_j)$, respectively calculated directly using one dimensional formulation of the scheme, e.g., when computing $A^{3\pm}_{x}$, we fix $y$ and apply 1D scheme in $x$-direction. For example, when $A_x$ is periodic in $x$-direction, we get
\begin{align*}
& A_{x}^{-} \approx \alpha_{1}\mathcal{D}_{L}^{p} [A(\cdot,y),\alpha_{1}](x) + \alpha_1 \sum_{p=2}^{k} \sigma_{i,L}^{p-2} \mathcal{D}_{L}^{p} [A(\cdot,y),\alpha_{1}](x),\\
& A_{x}^{+}\approx -\alpha_{1} \mathcal{D}_{R}[A(\cdot,y),\alpha_{1}](x) -\alpha_{1}\sum_{p=2}^{k} \sigma_{i,R}^{p-2} \mathcal{D}_{R}^{p} [A(\cdot,y),\alpha_{1}](x).
\end{align*}
Here, we choose $\alpha_{1}=\beta/(c_{1}\Delta t)$. Similarly, to approximate $A_y^\pm$, we fix $x$ and obtain
\begin{align*}
& A_{y}^{-} \approx \alpha_{2} \mathcal{D}_{L} [A(x,\cdot),\alpha_{2}](y) + \alpha_{2}\sum_{p=2}^{k} \sigma_{i,L}^{p-2} \mathcal{D}_{L}^{p} [A(x,\cdot),\alpha_{2}](y),\\
& A_{y}^{+}\approx  -\alpha_{2} \mathcal{D}_{R} [A(x,\cdot),\alpha_{2}](y) -\alpha_{2}\sum_{p=2}^{k} \sigma_{i,R}^{p-2}\mathcal{D}_{R}^{p} [A(x,\cdot),\alpha_{2}](y),
\end{align*}
with $\alpha_{2}=\beta/(c_{2}\Delta t)$.
Note that in the 2D case,  we need to choose $\beta_{max}$ as half of that for the 1D case to ensure the unconditional stability of the scheme.

If it is non periodic boundary condition, we still use extrapolation with suitable order of accuracy for the derivative values at an outflow boundary, as in the 1D formulation. For the details see \cite{christlieb2019kernel}.

\subsection{3D magnetic potential equation}

Although the evolution equation for the 3D magnetic potential \eqref{eq:3dmagneticP} is significantly different from the evolution equation for the 2D scalar magnetic potential \eqref{eq:2dmagneticP}, we can still directly apply the scheme presented in section \eqref{sec:1dformulation} to 3D case. In our previous paper \cite{christlieb2014finite}, we needed an artificial resistivity term in 3D case. However, with the kernel based approach we no longer need those artificial terms in our equations, since our method is fully implicit.

Writing out the magnetic vector potential equation derived in section \eqref{sec:3DMP} in component form we have:
\begin{equation}
\partial_t A^1 = u^2 (\partial_x A^2) + u^3 (\partial_x A^3) - u^2 (\partial_y A^1) - u^3 (\partial_z A^1)
\end{equation}
\begin{equation}
\partial_t A^2 = - u^1 (\partial_x A^2) + u^1 (\partial_y A^1) + u^3 (\partial_y A^3) - u^3 (\partial_z A^2)
\end{equation}
\begin{equation}
\partial_t A^3 = - u^1 (\partial_x A^3) - u^2 (\partial_y A^3) + u^1 (\partial_z A^1) + u^2 (\partial_z A^2)
\end{equation}
with $\Av = (A^1, A^2, A^3)$ and $\uv=(u^1,u^2,u^3)$.
While $\Av$ in 3D is not strictly a H-J equation, with the new implicit approach we can simply apply the ideas from the 2D case. Then we can obtain the following equations using the Lax-Friedrichs flux splitting:
\begin{subequations}
	\begin{align} 
	\frac{dA^{1}(t)}{dt} = & \quad u^2 \frac{(A^{2+}_{x}+A^{2-}_{x})}{2} -c_2 \frac{(A^{2-}_{x}-A^{2+}_{x})}{2}
	+ u^3 \frac{(A^{3+}_{x}+A^{3-}_{x})}{2} - c_3 \frac{(A^{3-}_{x}-A^{3+}_{x})}{2}  \notag\\
	& - u^2 \frac{(A^{1+}_{y}+A^{1-}_{y})}{2}-c_2 \frac{(A^{1-}_{y}-A^{1+}_{y})}{2}
	- u^3 \frac{(A^{1+}_{z}+A^{1-}_{z})}{2}-c_3 \frac{(A^{1-}_{z}-A^{1+}_{z})}{2},    \\
	\frac{dA^{2}(t)}{dt} =& -u^1 \frac{(A^{2+}_{x}+A^{2-}_{x})}{2}-c_1 \frac{(A^{2-}_{x}-A^{2+}_{x})}{2}
	+ u^1 \frac{(A^{1+}_{y}+A^{1-}_{y})}{2} -c_1 \frac{(A^{1-}_{y}-A^{1+}_{y})}{2} \notag\\
	& + u^3 \frac{(A^{3+}_{y}+A^{3-}_{y})}{2} -c_3 \frac{(A^{3-}_{y}-A^{3+}_{y})}{2}
	- u^3 \frac{(A^{2+}_{z}+A^{2-}_{z})}{2}-c_3 \frac{(A^{2-}_{z}-A^{2+}_{z})}{2}, \\
	\frac{dA^{3}(t)}{dt} = & - u^1 \frac{(A^{3+}_{x}+A^{3-}_{x})}{2}-c_1 \frac{(A^{3-}_{x}-A^{3+}_{x})}{2}
	- u^2 \frac{(A^{3+}_{y}+A^{3-}_{y})}{2}-c_2 \frac{(A^{3-}_{y}-A^{3+}_{y})}{2}  \notag\\
	& + u^1 \frac{(A^{1+}_{z}+A^{1-}_{z})}{2} -c_1 \frac{(A^{1-}_{z}-A^{1+}_{z})}{2} 
	+ u^2 \frac{(A^{2+}_{z}+A^{2-}_{z})}{2} -c_2 \frac{(A^{2-}_{z}-A^{2+}_{z})}{2} ,
	\end{align}  
\end{subequations}
at $(x_i, y_j, z_k)$,
where 
\begin{equation*}
	c_1 = \max_{i,j,k}{|u^1_{i,j,k}|},   \quad 
	c_2 = \max_{i,j,k}{|u^2_{i,j,k}|},  \quad  \text{and}  \quad  
	c_3 = \max_{i,j,k}{|u^3_{i,j,k}|}.
	\end{equation*}

Similarly, we use 1D 
kernel-based formulation to approximate the derivatives of the magnetic vector potential components  $A^{*\pm}_{x}$,  $A^{*\pm}_{y}$, and  $A^{*\pm}_{z}$, where $*$ denotes $1, 2$ and $3$.
In addition, since the density or pressure may become negative in some problems such as the blast wave problem, we use the positivity preserving limiter idea developed in our previous work \cite{christlieb2014finite}. For more details on positivity limiter, see \cite{christlieb2016high, christlieb2015positivity, seal2016explicit}
\section{Central finite difference discretization of $\nabla \times \Av$}

According to the Constrained Transport formulation, we applied a discrete curl operator to the magnetic potential at each stage of CT framework for the purpose of obtaining a divergence free magnetic field. Here we will present the strategy we used to approximate the curl operator.

\subsection{Curl in 2D}

We apply a discrete version of 2D curl to the equation \eqref{eq:2drelation} as 
	\begin{equation}
	B^{1}_{i,j} : = D_{y}|_{i,j}A^3 \quad\quad \text{and} \quad\quad  B^{2}_{i,j} : = - D_{x}|_{i,j}A^3,
	\end{equation} 
where $D_x$ and $D_y$ are the discrete versions of the partial derivatives $\partial_x$ and $\partial_y$. We use these version of curls to satisfy the discrete divergence free constraint as follows
	\begin{equation*}
	\nabla \cdot B_{i,j} : = D_{x}|_{i,j}B^1 + D_{y}|_{i,j}B^2 
	= D_{x}D_{y}|_{i,j}A^3 - D_{y}D_{x}|_{i,j}A^3 = 0.
	\end{equation*}

Here, we choose 4th-order central finite differences to get high order accuracy
\begin{subequations}
	\begin{align}
	& D_{x}|_{i,j}A^3 : = \frac{1}{12\Delta x}(A_{i-2, j}^3 - 8 A_{i-1, j}^3 + 8A_{i+1, j}^3 - A_{i+2, j}^3 ),\\
	& D_{y}|_{i,j}A^3 : = \frac{1}{12\Delta y}(A_{i, j-2}^3 - 8 A_{i, j-1}^3 + 8A_{i, j+1}^3 - A_{i, j+2}^3 ),
	\end{align}
\end{subequations}
where $A=A^3$ is only third component of magnetic potential.

\subsection{Curl in 3D}

We use the following discrete version of the 3D curl
\begin{subequations}
	\begin{align}
	& B^1_{i,j,k} := D_{y}|_{i,j,k}A^3 - D_{z}|_{i,j,k}A^2, \\
	& B^2_{i,j,k} := D_{z}|_{i,j,k}A^1 - D_{x}|_{i,j,k}A^3, \\
	& B^3_{i,j,k} := D_{x}|_{i,j,k}A^2 - D_{y}|_{i,j,k}A^1, 
	\end{align}
\end{subequations}
where $D_x$, $D_y$, and $D_z$ are notations for the discrete versions of partial derivatives $\partial_x$ , $\partial_y$ and $\partial_z$, respectively. With these versions of curls, we get the following divergence free condition satisfied
	\begin{equation}
	\begin{aligned}
	\nabla \cdot B_{i,j,k} :  
	& = D_{x}|_{i,j,k}B^1 + D_{y}|_{i,j,k}B^2 + D_{z}|_{i,j,k}B^3 \\
	& = D_{x} D_{y}|_{i,j,k}A^3 - D_{x}D_{z}|_{i,j,k}A^2 + D_{y}D_{z}|_{i,j,k}A^1 - D_{y}D_{x}|_{i,j,k}A^3 \\
	& + D_{z}D_{x}|_{i,j,k}A^2 - D_{z}D_{y}|_{i,j,k}A^1 \\
	& = 0
	\end{aligned}
	\end{equation}
As in 2D case, we use 4th-order central finite differences to get high order accuracy
\begin{subequations}
	\begin{align}
	& D_{x}|_{i,j,k}A^{*} : = \frac{1}{12\Delta x}(A^{*}_{i-2, j, k} - 8 A^{*}_{i-1, j, k} + 8A^{*}_{i+1, j, k} - A^{*}_{i+2, j, k}), \\
	& D_{y}|_{i,j,k}A^{*} : = \frac{1}{12\Delta y}(A^{*}_{i, j-2, k} - 8 A^{*}_{i, j-1, k} + 8A^{*}_{i, j+1, k} - A^{*}_{i, j+2, k}), \\
	& D_{z}|_{i,j,k}A^{*} : = \frac{1}{12\Delta z}(A^{*}_{i, j, k-2} - 8 A^{*}_{i, j, k-1} + 8A^{*}_{i, j, k+1} - A^{*}_{i, j, k+2}).
	\end{align}
\end{subequations}
\section{Numerical Results}

In this section, we present the numerical results to demonstrate the accuracy and efficiency of the new method. 
We use third order SSP RK method for time discretization. Time step is chosen as 
\begin{align}
\Delta t=\frac{\text{CFL}}{em*cd}
\end{align}
where, the CFL number is 0.5, $em = {\max(\lambda ^5, \lambda ^6)}$, $cd = \frac{1}{\Delta x} + \frac{1}{\Delta y}$ for 2D and $cd = \frac{1}{\Delta x} + \frac{1}{\Delta y} + \frac{1}{\Delta z}$ for 3D.

\subsection{Smooth vortex test in MHD}
\qquad We first test the smooth vortex problem in 2D with non zero magnetic field to show the accuracy of the method within the constrained transport formulation.The initial conditions are 
	\begin{equation*}
	(\rho, u^1, u^3, u^3, p, B^1, B^2, B^3) = (1, 1, 1, 0, 1, 0, 0, 0) 
	\end{equation*}
with perturbations on $u^1, u^2, B^1, B^2$ and $p$ as :
	\begin{align*}
	&(\delta u^1, \delta u^2) = \frac {\mu} {2\pi} e^{0.5(1-r^2)}(-y,x). \\
	&(\delta B^1, \delta B^2) = \frac {\kappa} {2\pi} e^{0.5(1-r^2)}(-y,x).  \\
	& \delta p = \frac {\mu ^y(1-r^2) - \kappa ^2} {8\pi ^2} e^{(1-r^2)}.
	\end{align*}
And the initial condition for magnetic potential is
	\begin{equation*}
	A^3(0, x, y) = \frac{\mu}{2\pi}e^{0.5(1-r^2)}
	\end{equation*}
where $r^2 = x^2 + y^2$. The vortex strength is taken as $\mu = 5.389489439$ and $\kappa = \sqrt{2}\mu$ such that the lowest pressure is around $5.3\times10^{-12}$ which happens in the center of the vortex. The domain is $[-10, 10]\times[-10,10]$ and periodic boundary condition is used on all four boundaries. In Table \ref{tab:ex1}, we present the errors of $\rho$ at $t=0.05$ with the mesh size $160\times160$, demonstrating that the scheme is third order as designed.


\begin{table}[htb]
	\caption{\label{tab:ex1} \em \large Smooth vortex problem. Errors of $\rho$ and orders of accuracy. } 
	\centering
	\vspace{0.1cm}   
	\begin{Large}    
		\begin{tabular}{|c|cc|cc|}
			\hline
		    $N_x\times N_y$&  $ \quad  L_1$ error & \quad  order  & $ \quad L_{\infty}$ error &  \quad order  \\\hline 
			$20\times20$  &  \quad 2.827E-03  &  \quad   --   & \quad 1.479E-01  &  \quad  --    \\
			$40\times40$  &  \quad 2.982E-04  &  \quad 3.245  & \quad 1.839E-02  &  \quad 3.007  \\
			$80\times80$  & \quad 1.861E-05  &  \quad 4.003  & \quad 1.327E-03  &  \quad 3.793  \\ 
			$160\times160$  & \quad 1.102E-06  & \quad 4.078  & \quad 1.126E-04  & \quad 3.559  \\
			$320\times320$  & \quad 7.260E-08  &  \quad 3.924  & \quad 1.170E-05 &  \quad 3.267  \\\hline	  
		\end{tabular}
	\end{Large}
\end{table}


\subsection{2D Orszag-Tang Vortex}
\begin{figure}[h]     
	\begin{center}
		\subfigure[Kernal-based method.] {\includegraphics[width=0.45\textwidth]{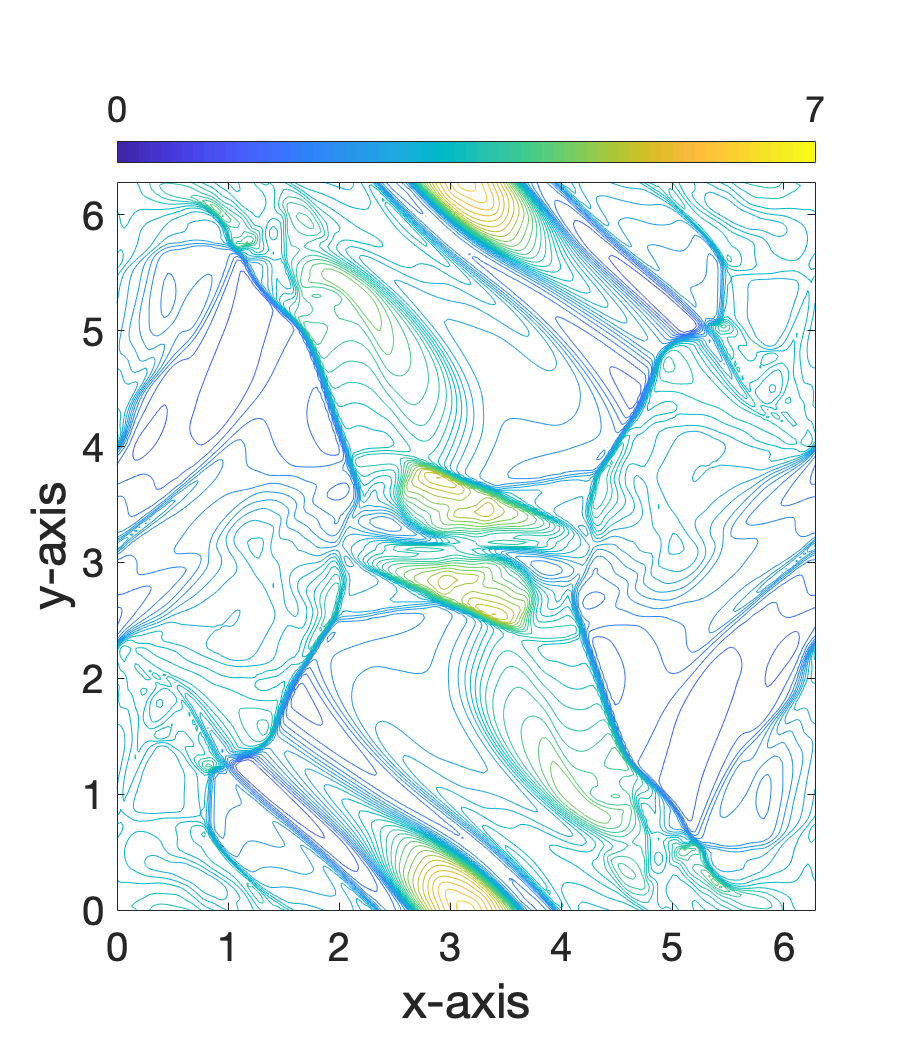}}
		\subfigure[Previous method \cite{christlieb2014finite}.] {\includegraphics[width=0.45\textwidth]{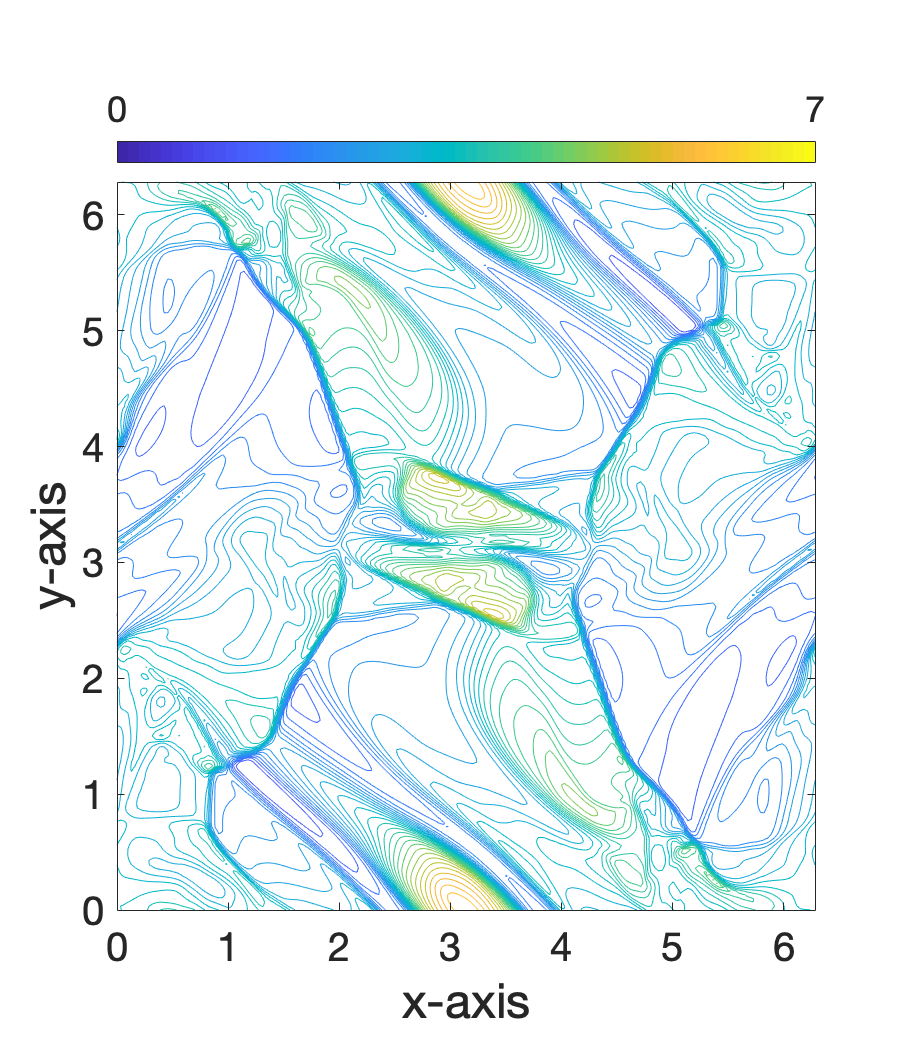}}
		\caption{\label{orz} \em Orszag-Tang vortex problem. Contour plots of density at $t = 3$ with $192\times192$ grid points.} 
	\end{center}
\end{figure}

 The Orszang-Tang vortex problem is a standard model problem for testing $\nabla \cdot \Bv = 0$ condition, since the solution of the problem is sensitive to divergence errors at late times. The initial conditions are given as
	\begin{equation*}
	(\rho, u^1, u^2, u^3, p, B^1, B^2, B^3) = (\gamma^2, -\sin(y), \sin(x), 0, \gamma, -\sin(y), \sin(2x),0 ), 
	\end{equation*}
where $\gamma = 5/3$ is the ideal gas constant and the initial magnetic potential:
	\begin{equation*}
	A^3 = 0.5\cos(2x) + \cos(y)
	\end{equation*}
The computational domain is $[0, 2\pi]\times[0, 2\pi]$ and periodic boundary conditions are used everywhere. We test the schemes with $192\times192$ grid points. Although the problem has smooth initial condition, as solution progresses, several shock waves and a vortex shape appear in the center of the domain. In Figure \ref{orz}, we show density $\rho$ at time $t=3$, and compare the results with our previous method \cite{christlieb2014finite}. We can see that they are in good agreement. 	

\subsection{Cloud Shock}
\begin{figure}[h]
	\begin{center}
		\subfigure[Kernal-based method.] {\includegraphics[width=0.45\textwidth]{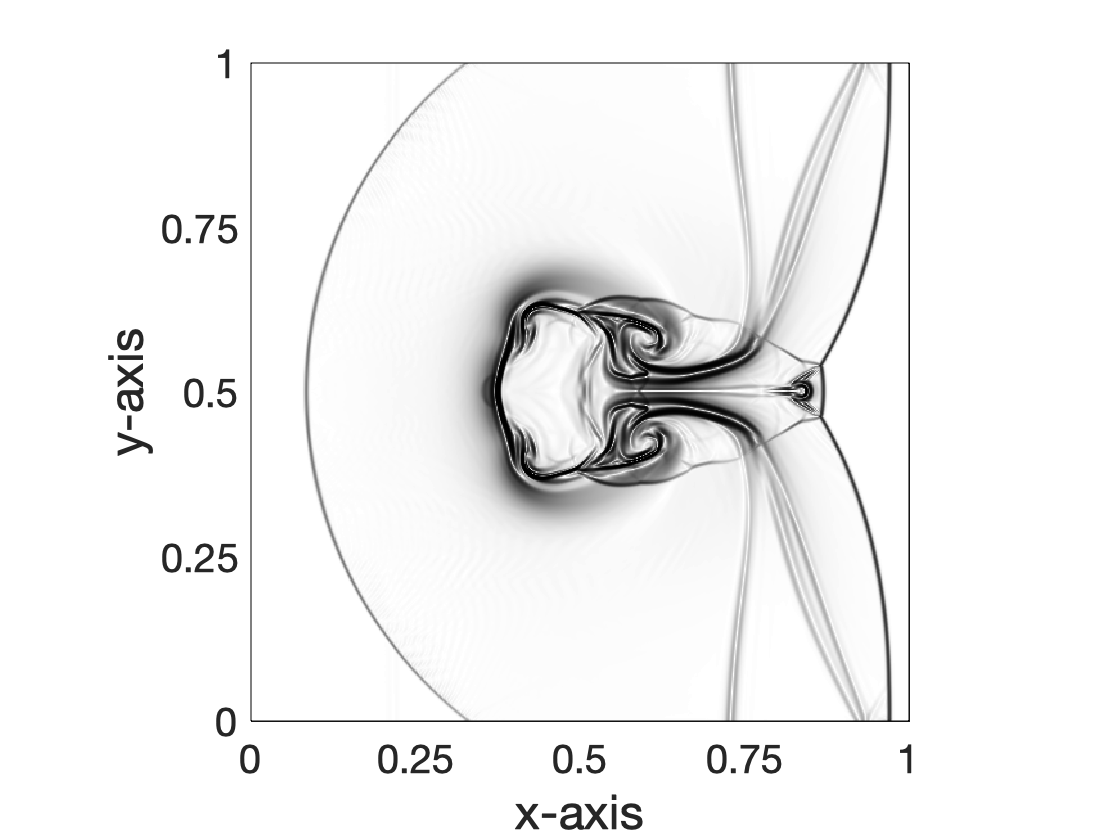}}
		\subfigure[Previous method \cite{christlieb2014finite}.] {\includegraphics[width=0.45\textwidth]{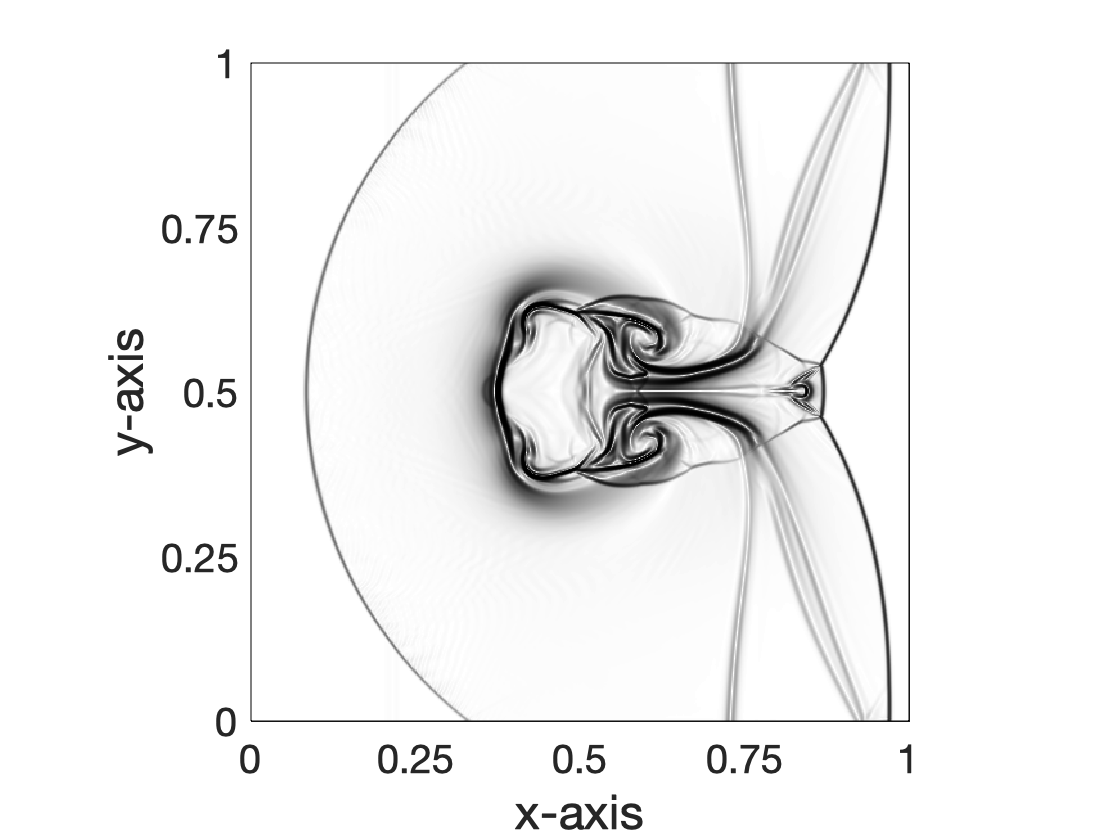}}
		\caption{\label{cloud} \em Cloud shock problem. Contour plots of $\| \Bv \|$ at $t = 0.06$ with $512\times512$ grid points.}
	\end{center}
\end{figure}

In this section we consider the 2D cloud-shock interaction problem, which models a strong shock passing  through a dense stationary bubble. The initial conditions include
	\begin{align*}
	 (\rho, u^1, u^2, u^3, p, B^1, B^2, B^3)  = 
	& \left\{\begin{array}{ll}
	(3.86859, 11.2536, 0, 0, 167.345, 0, 2.1826182, -2.1826182) & x < 0.05, \\
	\displaystyle (1, 0, 0, 0, 1, 0, 0.56418958, 0.56418958)  & x > 0.05, \\
	\end{array}
	\right.
	\end{align*}
and a circular cloud of density $\rho = 10$ and radius $r = 0.15$ centered at $(x, y) = (0.25, 0.5)$. The computational domain is $[0,1]\times[0,1]$ with mesh $512\times512$. We use inflow boundary condition at left boundary and outflow boundary condition elsewhere. 
The initial condition for magnetic potential:
\begin{align*}
	 A^3  = 
	& \left\{\begin{array}{ll}
	-2.1826182(x - 0.05) & x \leq 0.05, \\
	\displaystyle -0.56418958(x - 0.005)  & x \geq 0.05, \\
	\end{array}
	\right.
	\end{align*}
Figure \ref{cloud} presents Schlieren plots of $||\Bv||$ at $t=0.06$. The new method matches well with our numerical results of our previous method  \cite{christlieb2014finite}.

\subsection{2D Blast wave}
\begin{figure}[!h]
	\begin{center}
		\subfigure[density $\rho$] {\includegraphics[width=0.45\textwidth]{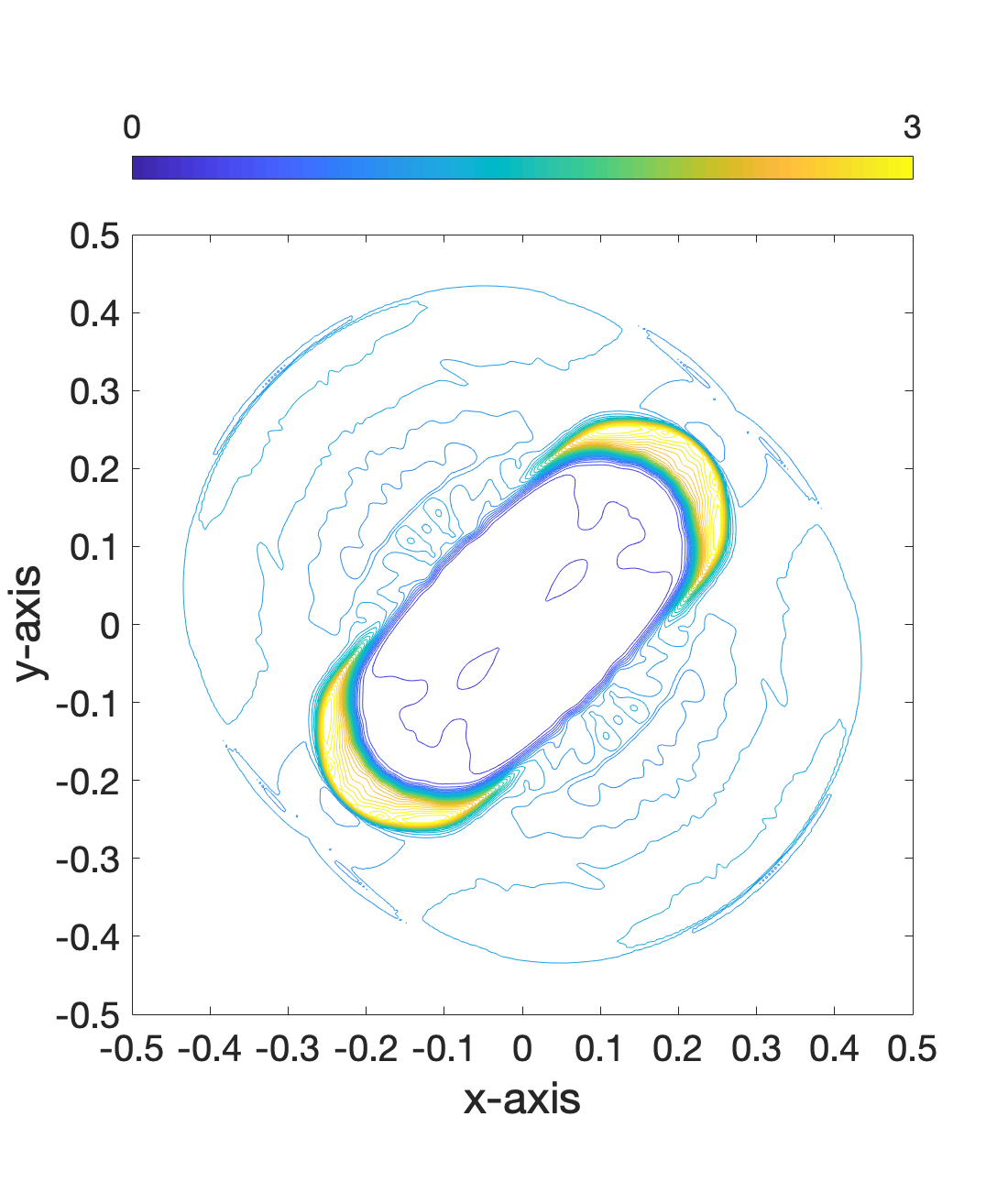}}
		\subfigure[pressure $p$] {\includegraphics[width=0.45\textwidth]{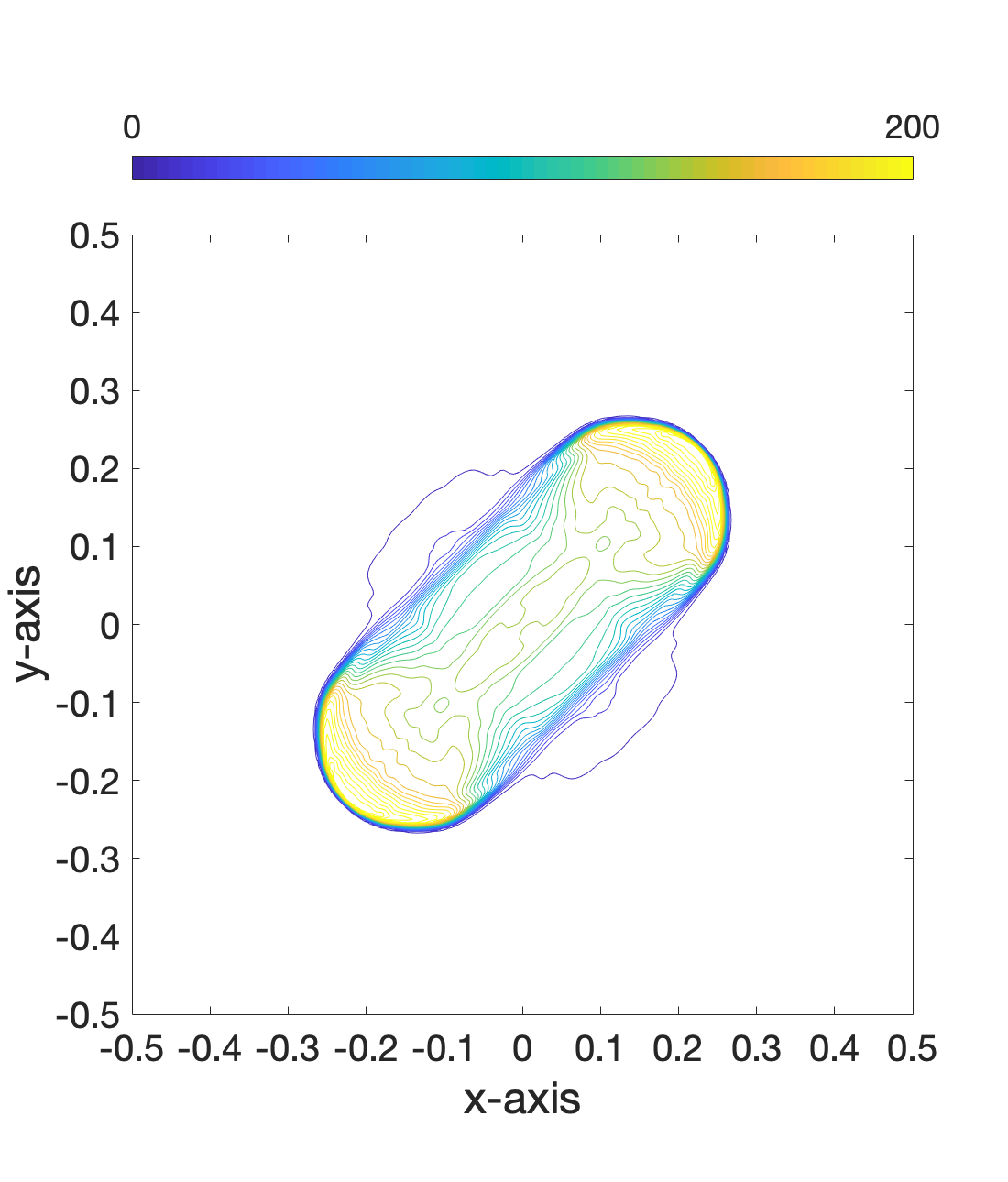}}
		\subfigure[$\| \uv \|$] {\includegraphics[width=0.45\textwidth]{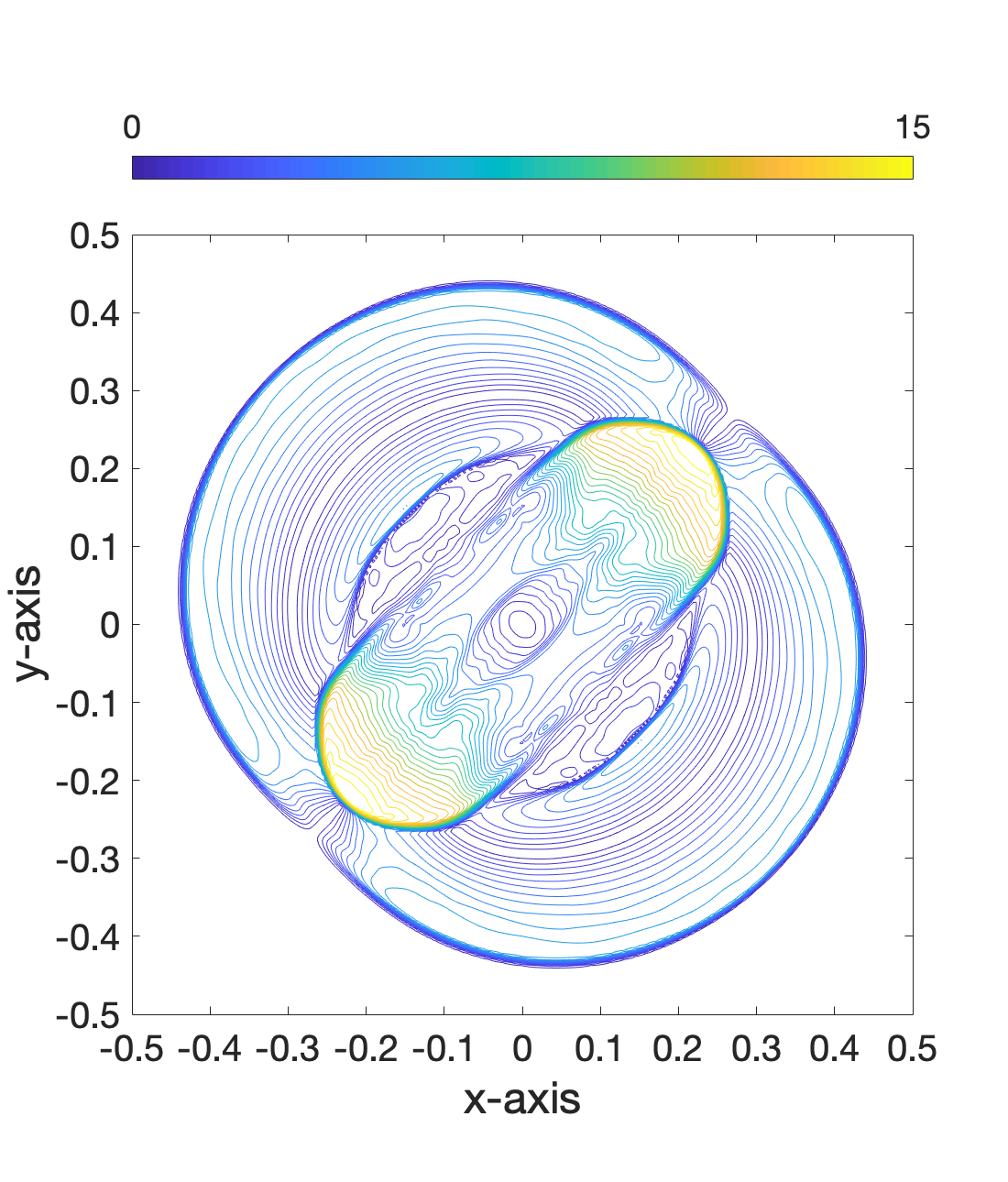}}
		\subfigure[$\| \Bv \|$] {\includegraphics[width=0.45\textwidth]{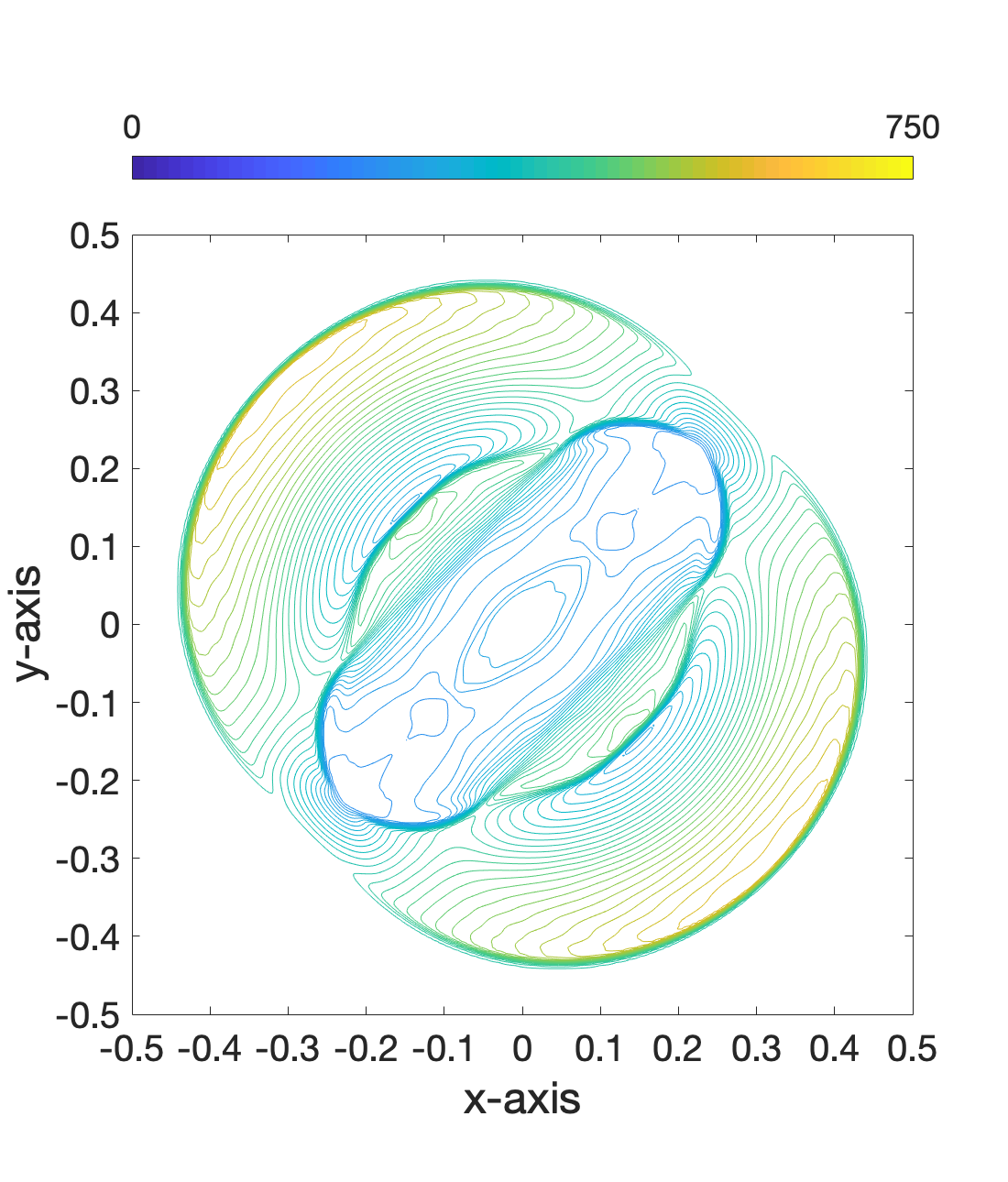}}
		\caption{\label{2dblast} \em 2D Blast wave problem. Contour plots at time $t = 0.01$ with $256\times256$ grid points. (a) density; (b) pressure; (c) the norm of $\uv$; (d) magnetic pressure $\|\Bv\|$.}
	\end{center}
\end{figure}

Here we investigate the 2D blast wave test problem, which has a strong shock causing negative density or pressure in simulations. To avoid negative density or pressure, we use the positivity preserving limiter we developed in reference \cite{christlieb2014finite}. The initial conditions are
	\begin{equation*}
	(\rho, u^1, u^2, u^3, B^1, B^2, B^3) = (1, 0, 0, 0, 50/{\sqrt{2\pi}}, 50/{\sqrt{2\pi}}, 0) 
	\end{equation*}
with a spherical pressure pulse
	\begin{align*}
	p =
	& \left\{\begin{array}{ll}
	1000 & r\leq 0.1, \\
	\displaystyle 0.1  & \text{otherwise}. \\
	\end{array}
	\right.
	\end{align*}	
The initial condition for magnetic potential can be obtained as
	\begin{equation*}
	A^3 = (0, 0,  50/{\sqrt{2\pi}} (y-x) ). 
	\end{equation*}
We use a domain $[-0.5, 0.5]\times[-0.5, 0.5]$ with $256\times256$ mesh. Outflow boundary conditions are applied everywhere. The results are shown in Figure \ref{2dblast} and they have good agreement with the results that we got in our previous paper. 

\subsection{2D Field Loop}
\begin{figure}[h]
	\begin{center}
		\subfigure[$\| \Bv \|$] {\includegraphics[width=0.45\textwidth]{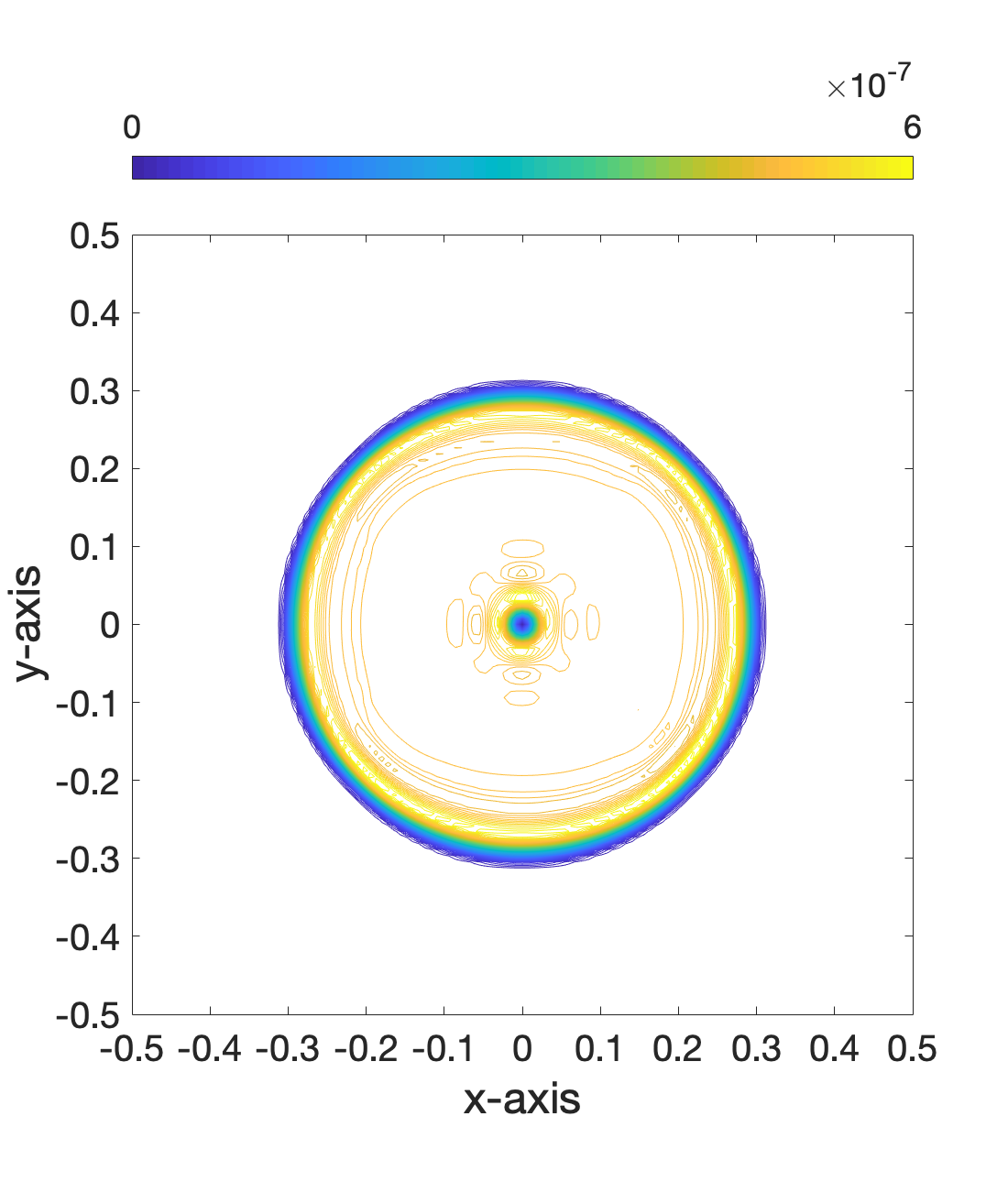}}
		\subfigure[$A$] {\includegraphics[width=0.45\textwidth]{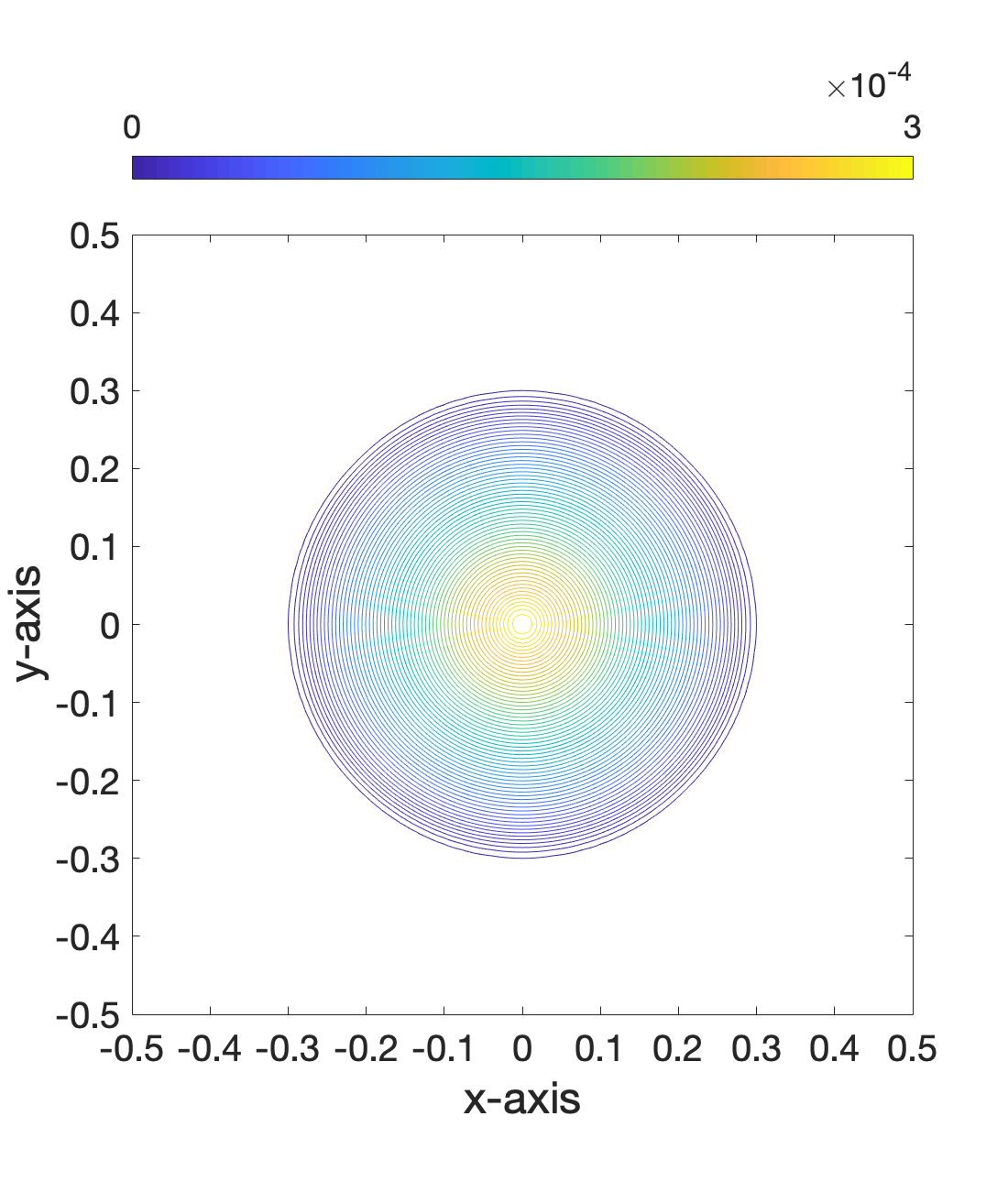}}
		\caption{\label{2dfield} \em 2D Field loop. Contour plots at time $t = 1$ with $128\times128$ grid points. a) Magnetic Pressure $\|\Bv\|$; b) Magnetic Potential $A^3$.}
	\end{center}
\end{figure}

The field loop problem is a strenuous test problem that moves a steady state. Here we test our recent method on this strenuous test case. The initial conditions are
	\begin{equation*}
	(\rho, u^1, u^2, p) = (1, \sqrt{5}\cos(\theta), \sqrt{5}\sin(\theta), 1) 
	\end{equation*}
with the angle $\theta = \arctan(0.5)$. The initial conditions for magnetic field $\Bv$ are determined by taking the curl of the magnetic potential, which is initialized with
	\begin{align*}
	A^3 =
	& \left\{\begin{array}{ll}
	0.001(R-r), & r\leq R, \\
	\displaystyle 0,  & \text{otherwise}, \\
	\end{array}
	\right.
	\end{align*}
where $r=\sqrt{x^2+y^2}$ and $R=0.3$. We use a domain $[-0.5, 0.5]\times[-0.5, 0.5]$ with $128\times128$ mesh. Periodic boundary conditions are applied to all sides. The loop is started in the middle of the domain and advected around the grid once, returned to initial location at time $t=1$ as shown in Figure \ref{2dfield}. We observe that the solution maintains the circular symmetry of the initial condition. 

\subsection{3D Field Loop}
\begin{figure}[h]
	\begin{center}
 		\subfigure[$\| \Bv \|$] {\includegraphics[width=0.45\textwidth]{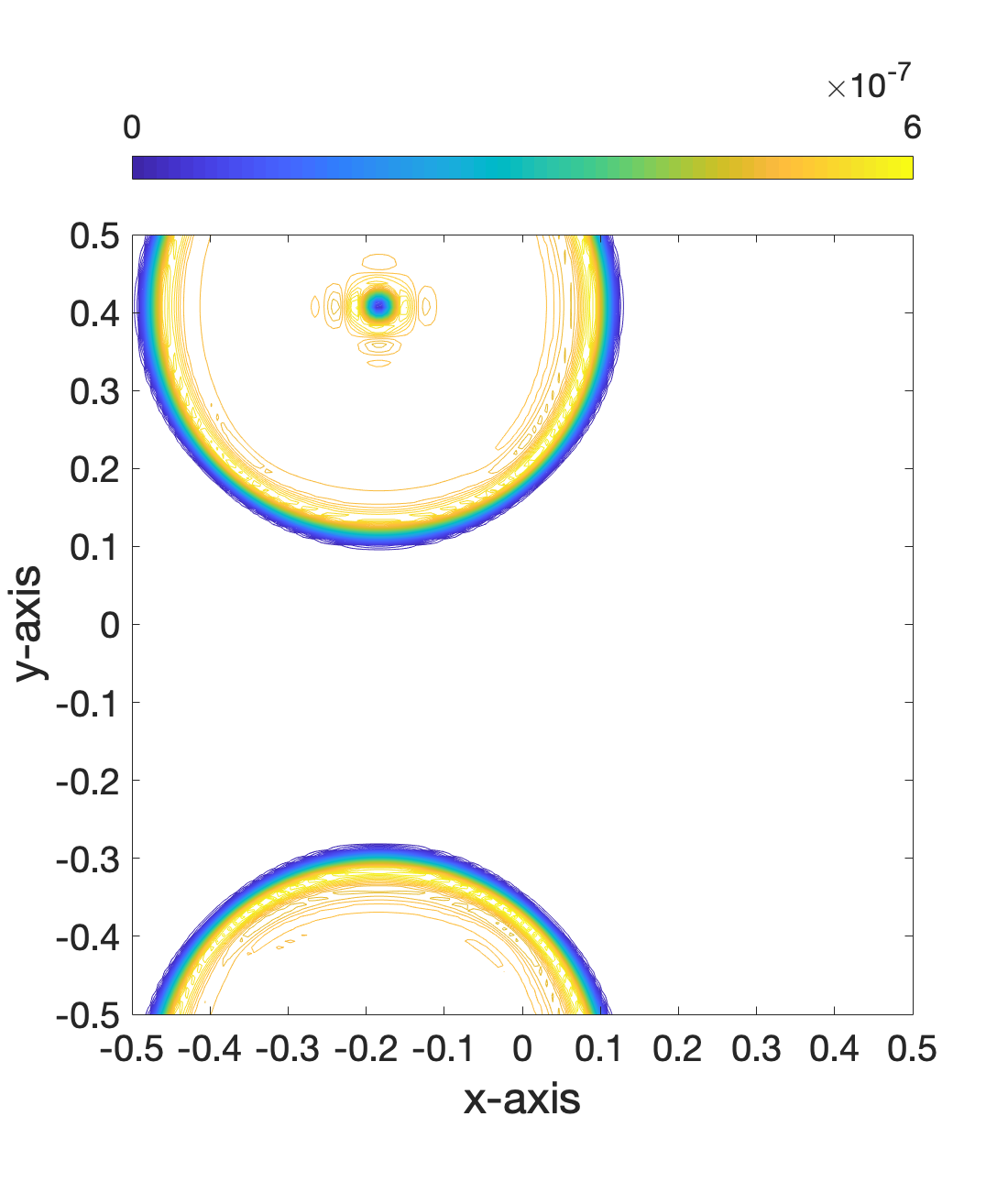}}
		\subfigure[$\Av$] {\includegraphics[width=0.45\textwidth]{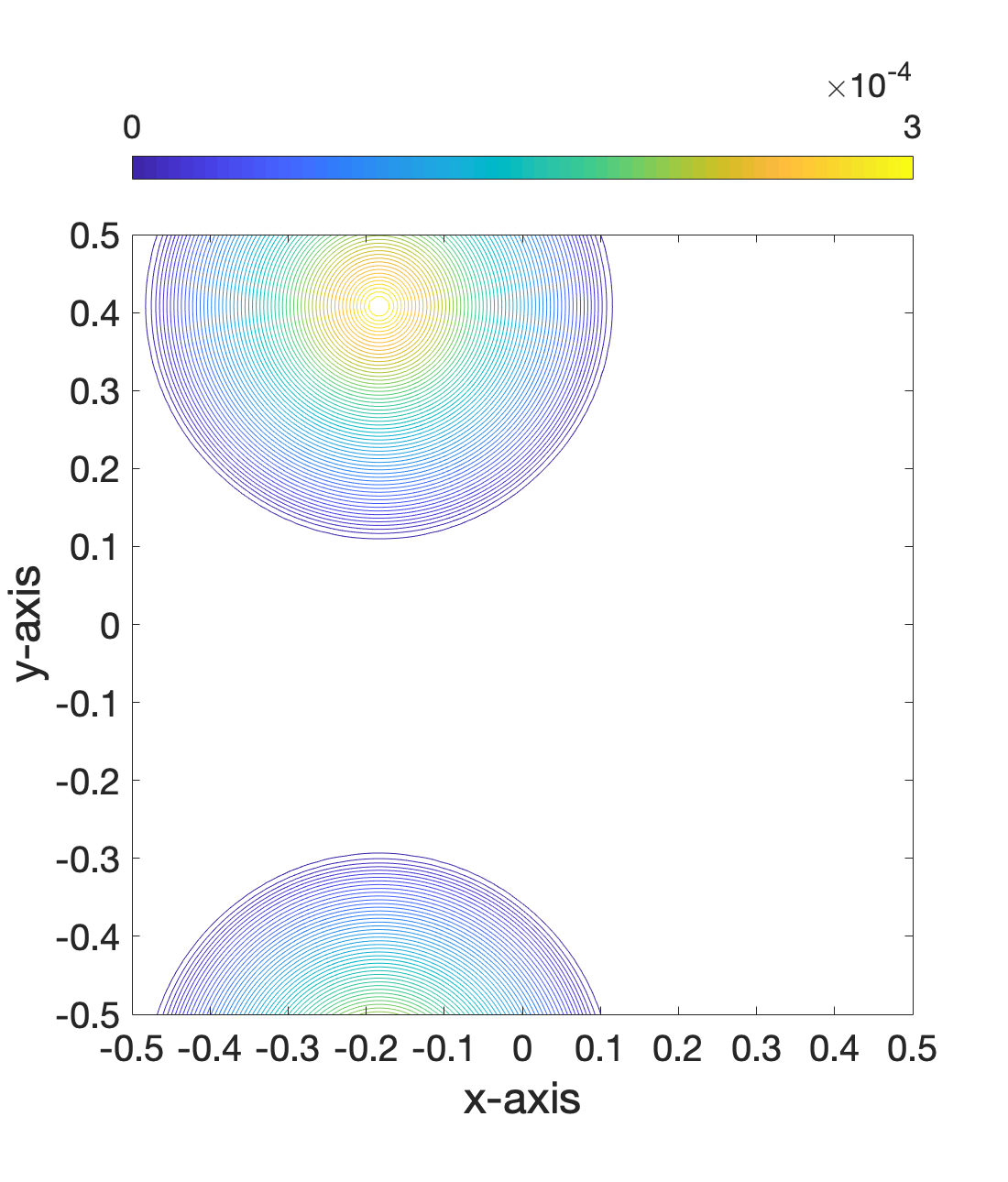}}
			\caption{\label{3dfield} \em 3D Field loop. Contour plots at time $t = 1$ with $128\times128$ grid points. The loop has been advected around the grid once. a) Magnetic Pressure; b) Magnetic Potential.}
		\end{center}
	\end{figure}

We tested an advecting field loop which moved diagonally across the boundary with an arbitrary initial angle.
The initial conditions are
	\begin{align*}
	(\rho, u^1, u^2, u^3, p) = (1, 2/{\sqrt{6}}, 1/{\sqrt{6}}, 1/{\sqrt{6}}, 1) 
	\end{align*}
The initial conditions for magnetic field are determined by taking the curl of the magnetic potential, which is given as magnetic potential :
	\begin{align*}
	A^3 =
	& \left\{\begin{array}{ll}
	0.001(R-r) & r\leq R \\
	\displaystyle 0  & \text{otherwise}, \\
	\end{array}
	\right.
	\end{align*}
where $A^1=0$, $A^2=0$, and $r=\sqrt{x^2+y^2}$ and $R=0.3$. 
We use a domain size $[-0.5, 0.5]\times[-0.5, 0.5]\times[-0.5, 0.5]$ with $128\times128\times128$ mesh. Periodic boundary conditions are applied to all sides. We observe that the field loop integrity is maintained, after advecting diagonally around the domain, until the final time, $t=1$. 
The results shown in Figure \ref{3dfield} are a 2D slice of the 3D solution taken at $z = 0$.

\subsection{3D Blast wave}
\begin{figure}[!ht]     
	\begin{center}
		\subfigure[density with kernel-based method] {\includegraphics[width=0.45\textwidth]{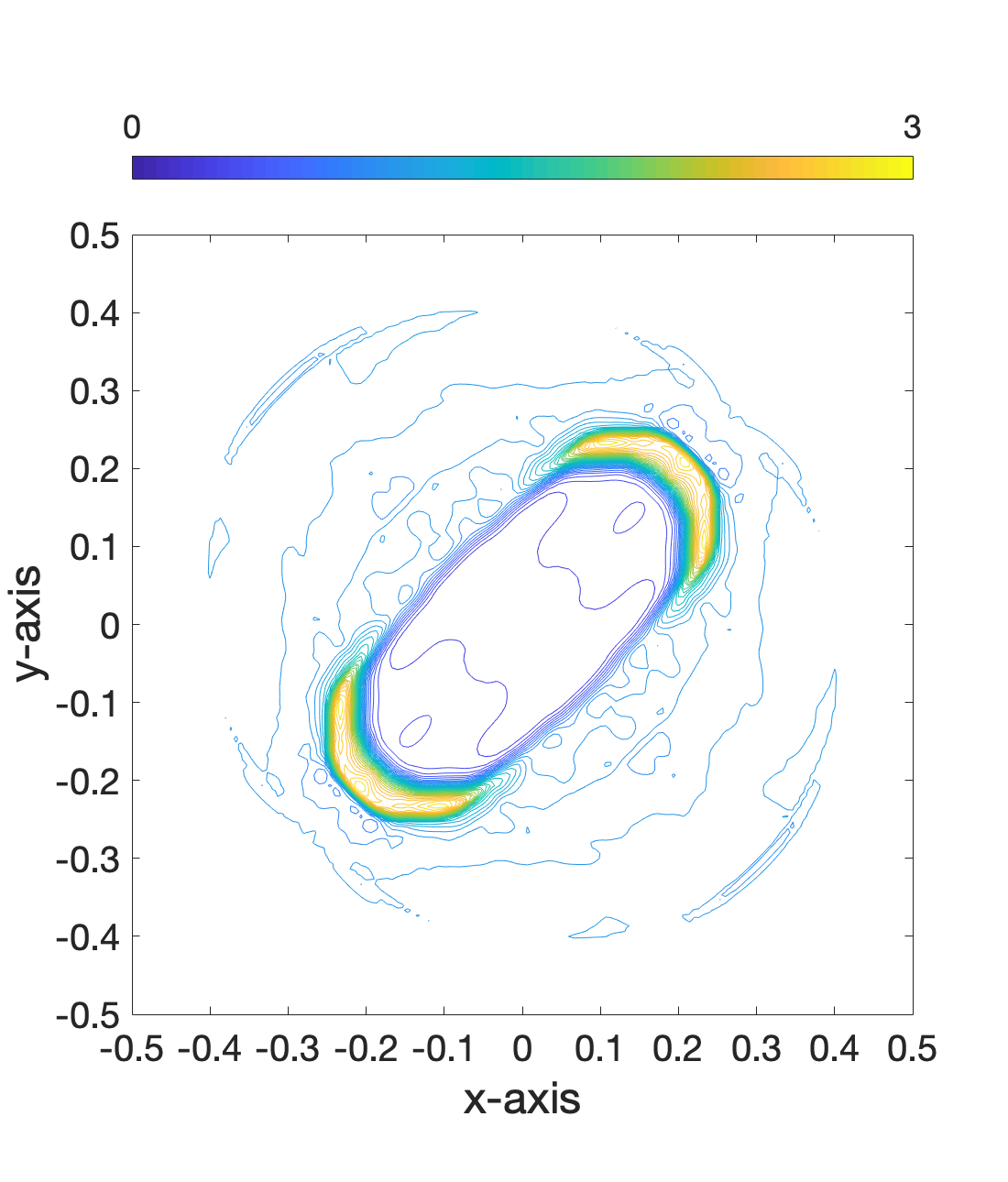}}
		\subfigure[density with previous method \cite{christlieb2014finite}.] {\includegraphics[width=0.45\textwidth]{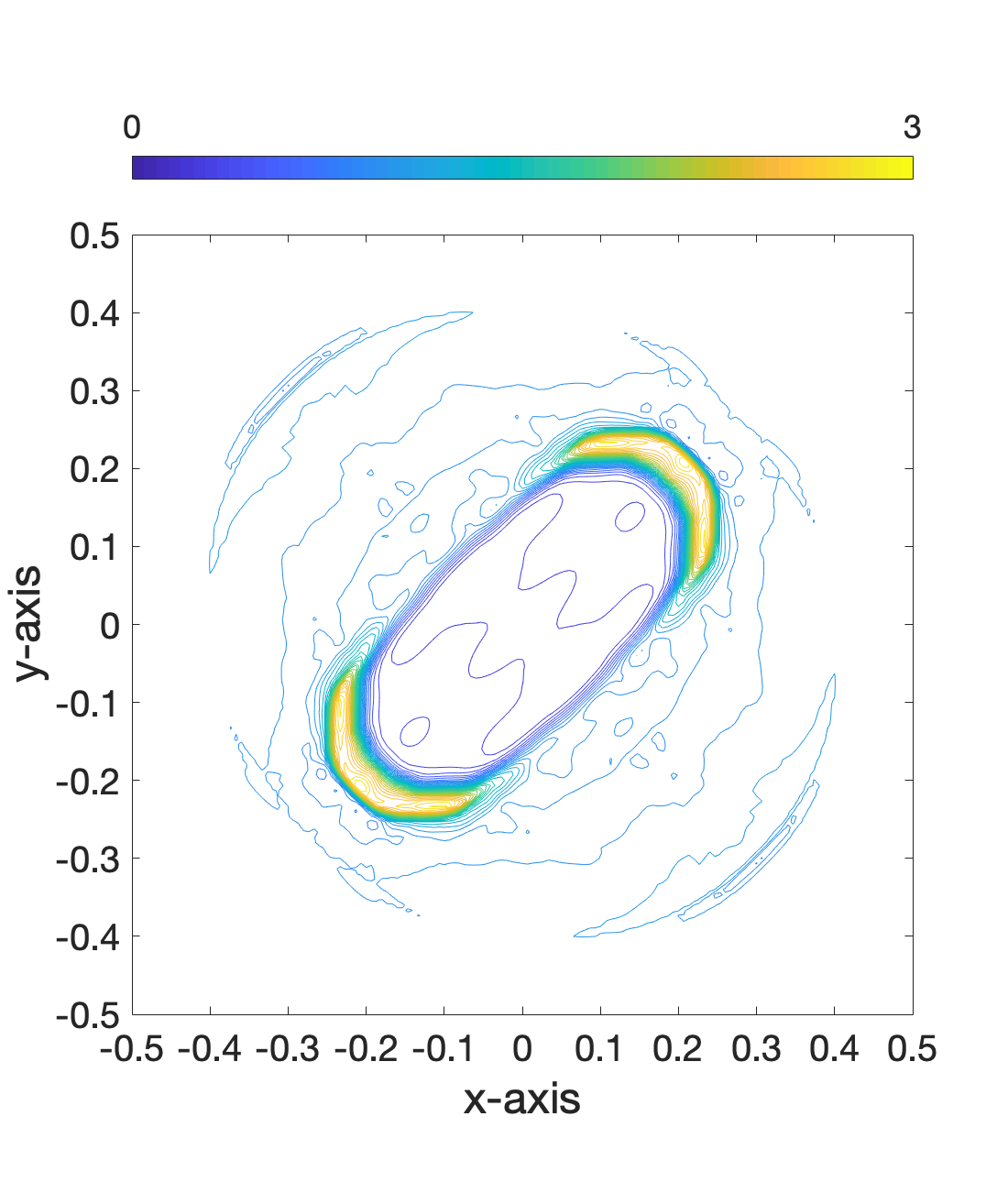}}
		\subfigure[pressure with kernel-based method] {\includegraphics[width=0.45\textwidth]{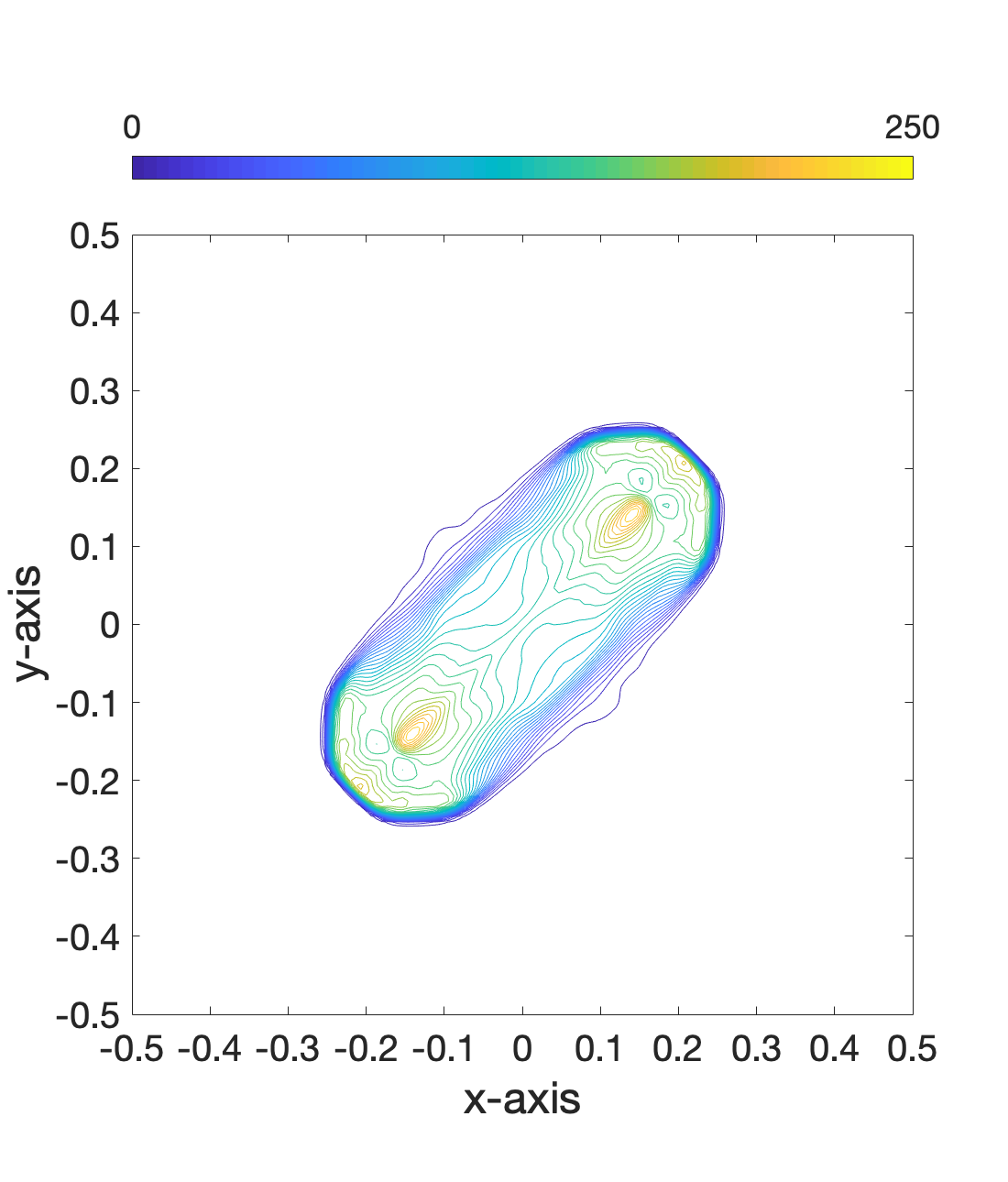}}
		\subfigure[pressure with previous method \cite{christlieb2014finite}] {\includegraphics[width=0.45\textwidth]{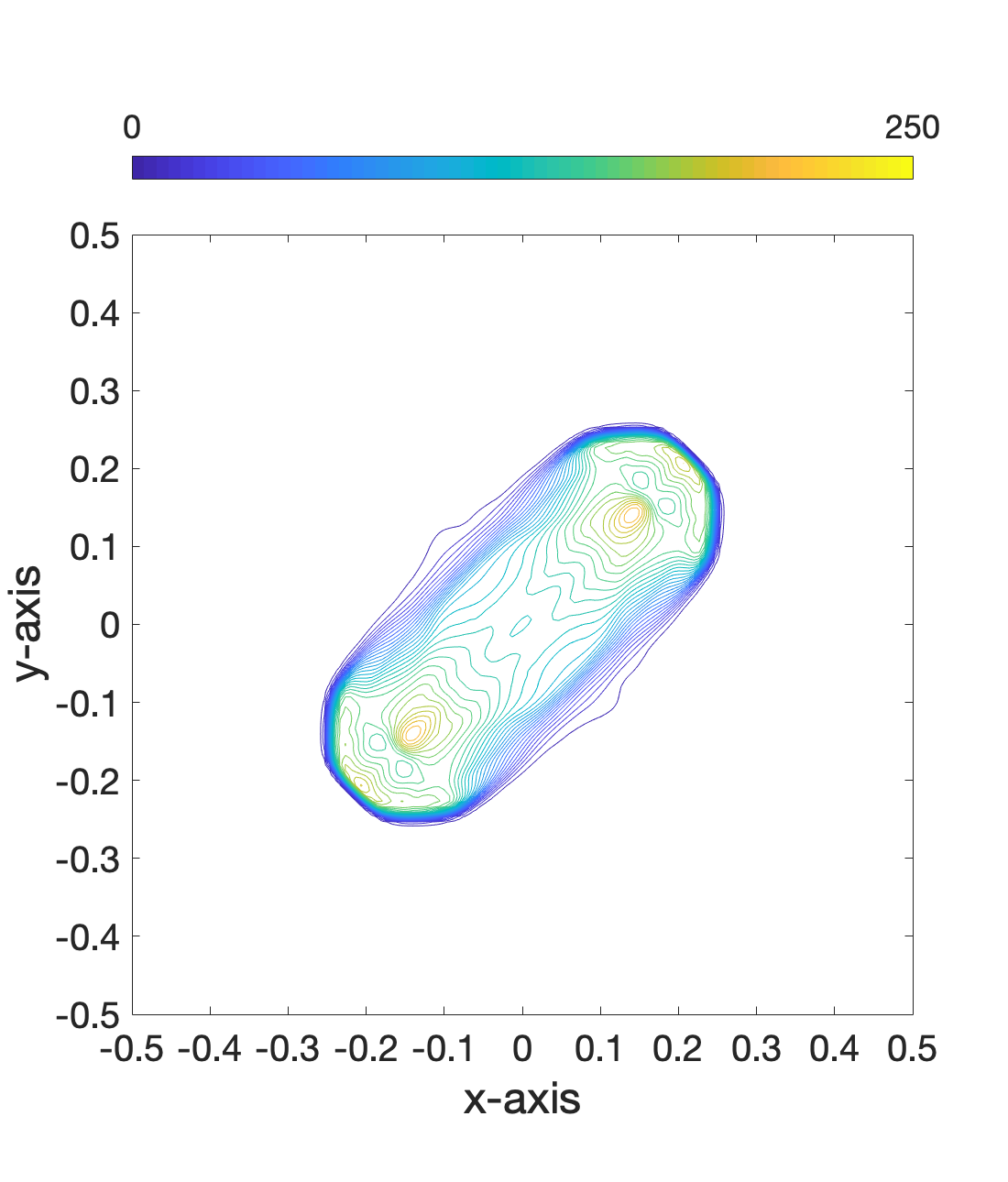}}
	\caption{\label{3dblast-1} \em 3D Blast wave problem. Contour plots at time $t = 0.01$ with $150\times150$ grid points}
	\end{center}
\end{figure}

\begin{figure}
	\begin{center}
		\subfigure[$\| \uv \|$ with kernel-based method] {\includegraphics[width=0.45\textwidth]{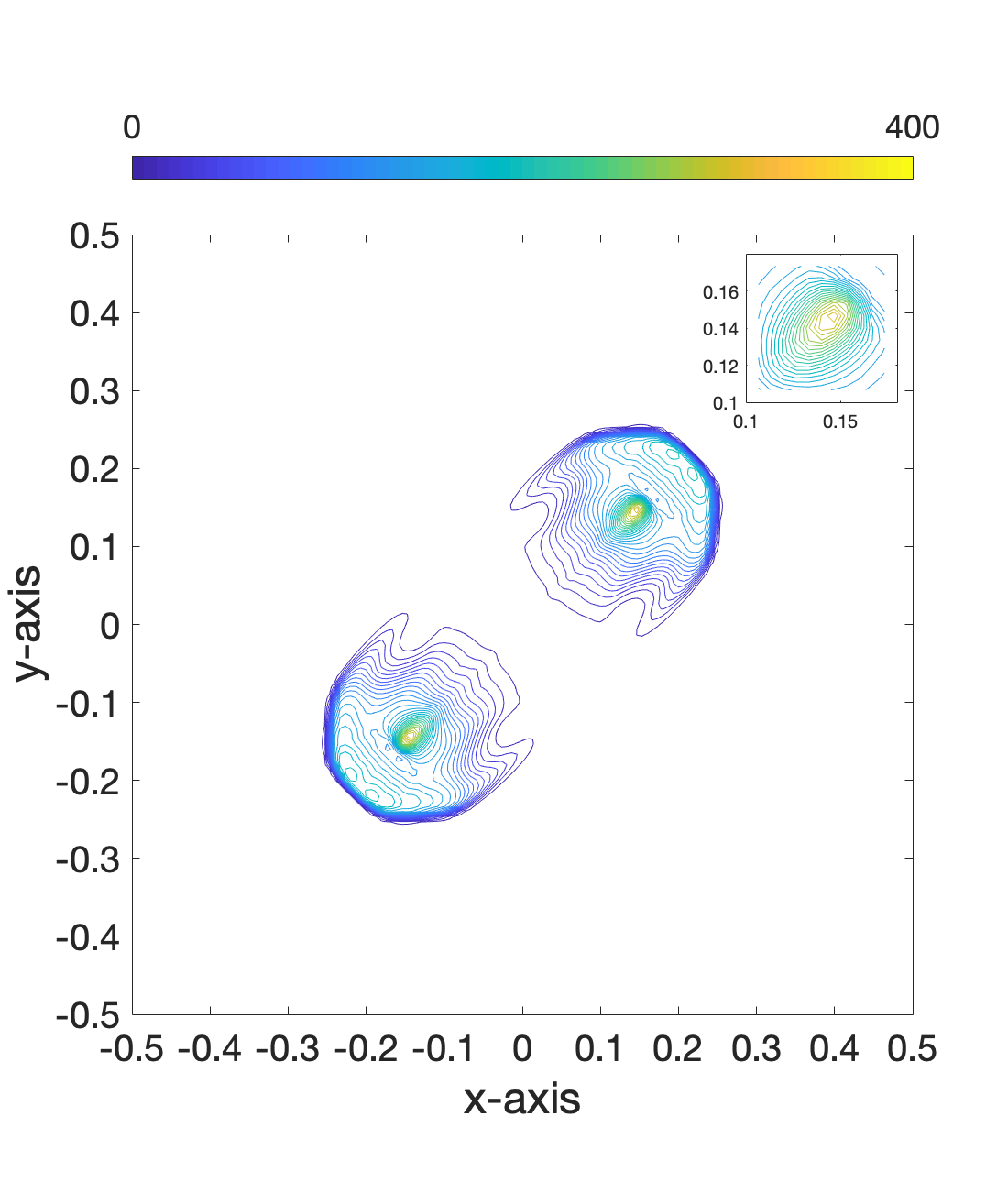}}
		\subfigure[$\| \uv \|$ with previous method \cite{christlieb2014finite}] {\includegraphics[width=0.45\textwidth]{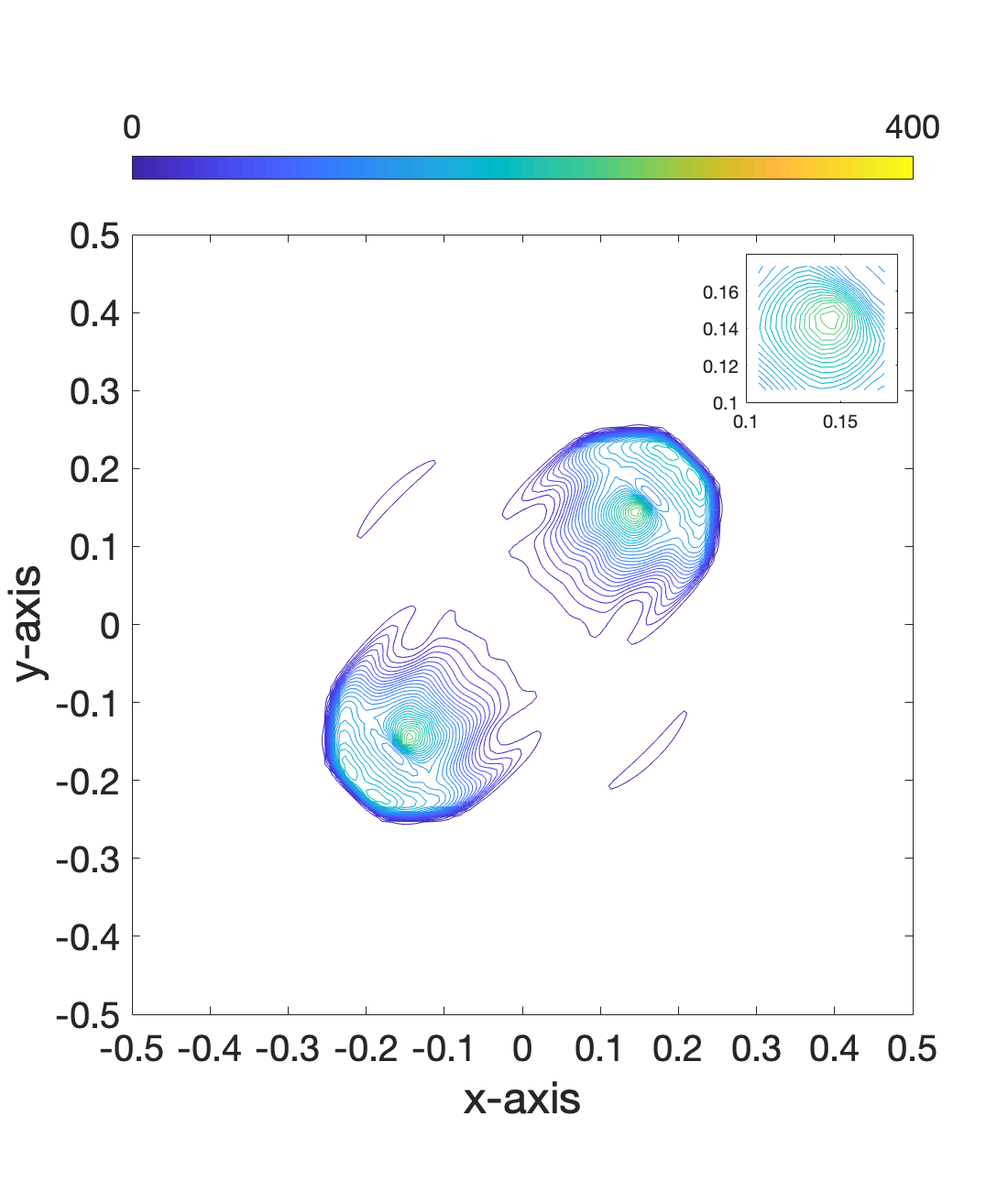}}
		\subfigure[$\| \Bv \|$ with kernel-based method] {\includegraphics[width=0.45\textwidth]{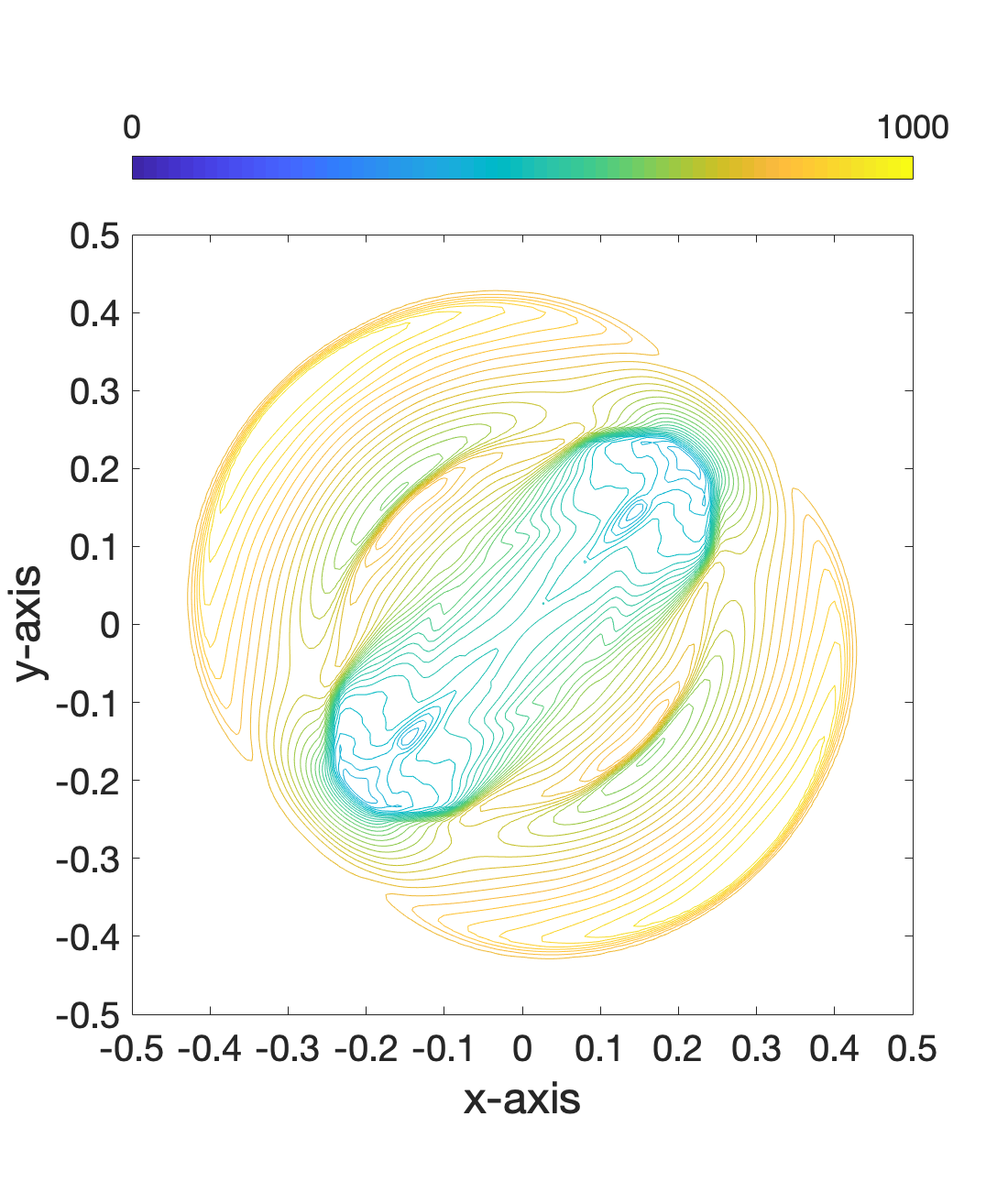}}
		\subfigure[$\| \Bv \|$ with previous method \cite{christlieb2014finite}] {\includegraphics[width=0.45\textwidth]{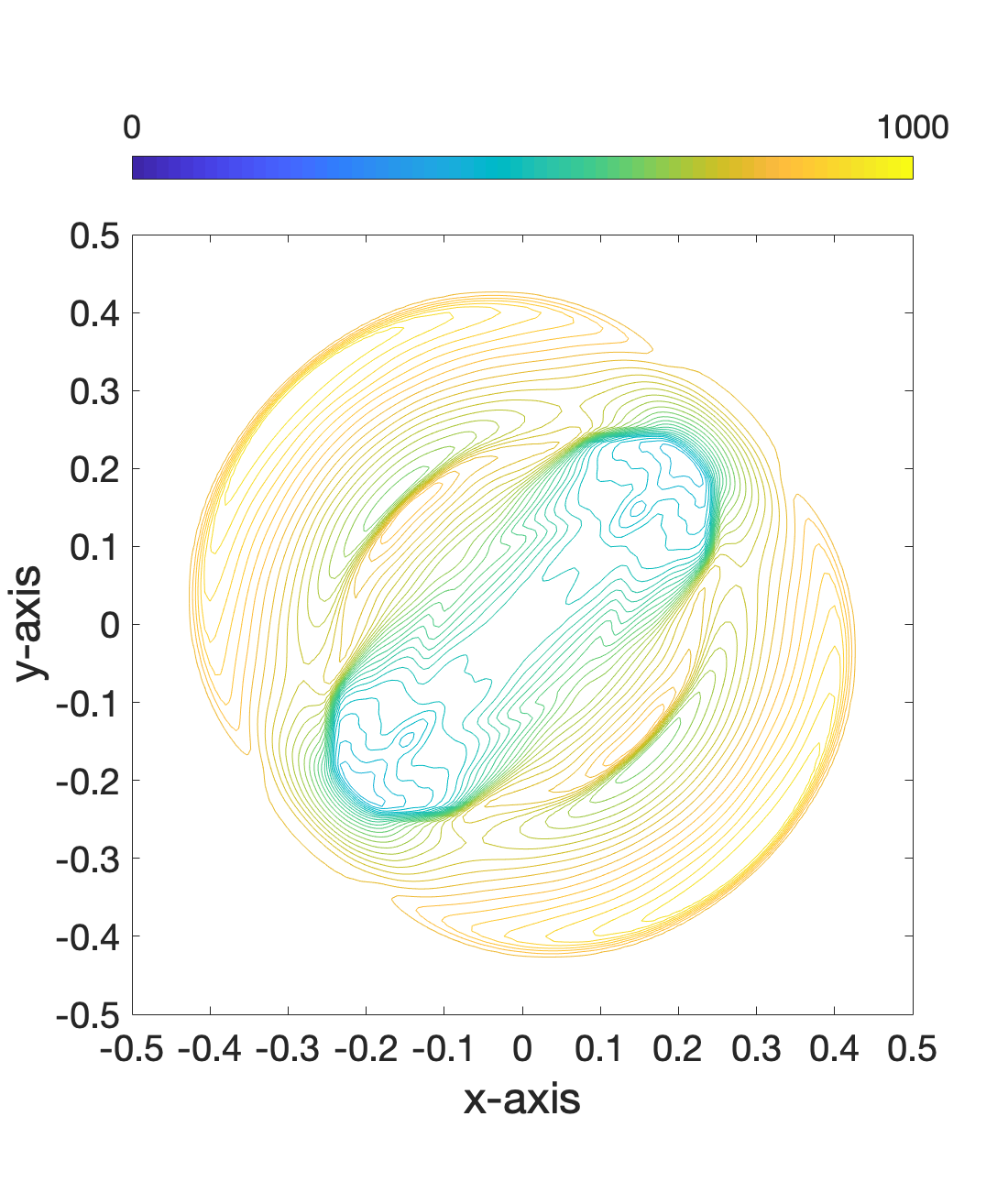}}
	\caption{\label{3dblast-2} \em 3D Blast wave problem. Contour plots at time $t = 0.01$ with $150\times150$ grid points. While the max $\| \uv \|$ value for previous method is 261, it is 320 for the kernel-based method. }
	\end{center}
\end{figure}

In this section we investigate the 3D version of the blast wave problem to show the strength of the new method, which eliminating the need for the diffusion limiter. The initial conditions are 
	\begin{equation*}
	(\rho, u^1, u^2, u^3, B^1, B^2, B^3) = (1, 0, 0, 0, 50/{\sqrt{2\pi}}, 50/{\sqrt{2\pi}}, 0) 
	\end{equation*}
with a spherical pressure pulse
	\begin{align*}
	p =
	& \left\{\begin{array}{ll}
	1000 & r\leq 0.1 \\
	\displaystyle 0.1  & \text{otherwise}. \\
	\end{array}
	\right.
	\end{align*}	
where $r=\sqrt{x^2+y^2+z^2}$. The initial condition for magnetic potential :
	\begin{equation*}
	\Av(0, x, y, z) = (0, 0,  50/{\sqrt{2\pi}}y - 50/{\sqrt{2\pi}}x). 
	\end{equation*}
We use a domain size $[-0.5, 0.5]\times[-0.5, 0.5]\times[-0.5, 0.5]$ with $150\times150\times150$ mesh. Outflow boundary conditions are applied everywhere.The results presented in Figure \ref{3dblast-1} and Figure \ref{3dblast-2} are the solutions cut at $z = 0$.
The sub plots highlight the differences between the two methods.  The key differences are the new method has a maximum value that achieves a peak of 320, which is 23\% higher than the old method, and is far less isotropic around the peek than the old method.  As demonstrated in a latter test, the isotropic behavior and lower the peek in the strong bast wave is due to the diffusion limiter that was needed with the old explicit constrained transport method.

	
\subsection{Kernel-based method with Diffusion Terms}

\begin{figure}[!hbt]    
	\begin{center}
		\subfigure[ density ] {\includegraphics[width=0.45\textwidth]{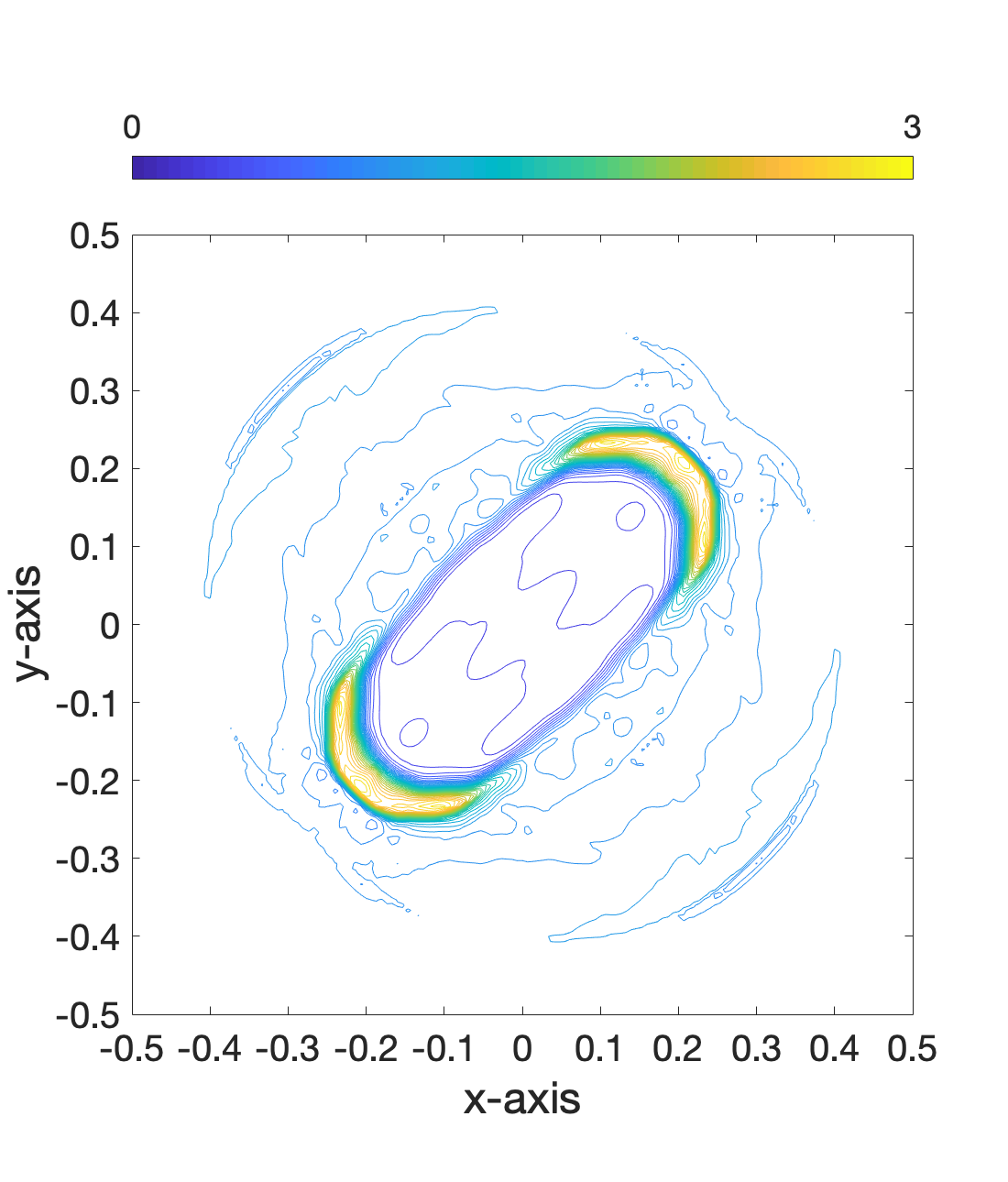}}
		\subfigure[pressure] {\includegraphics[width=0.45\textwidth]{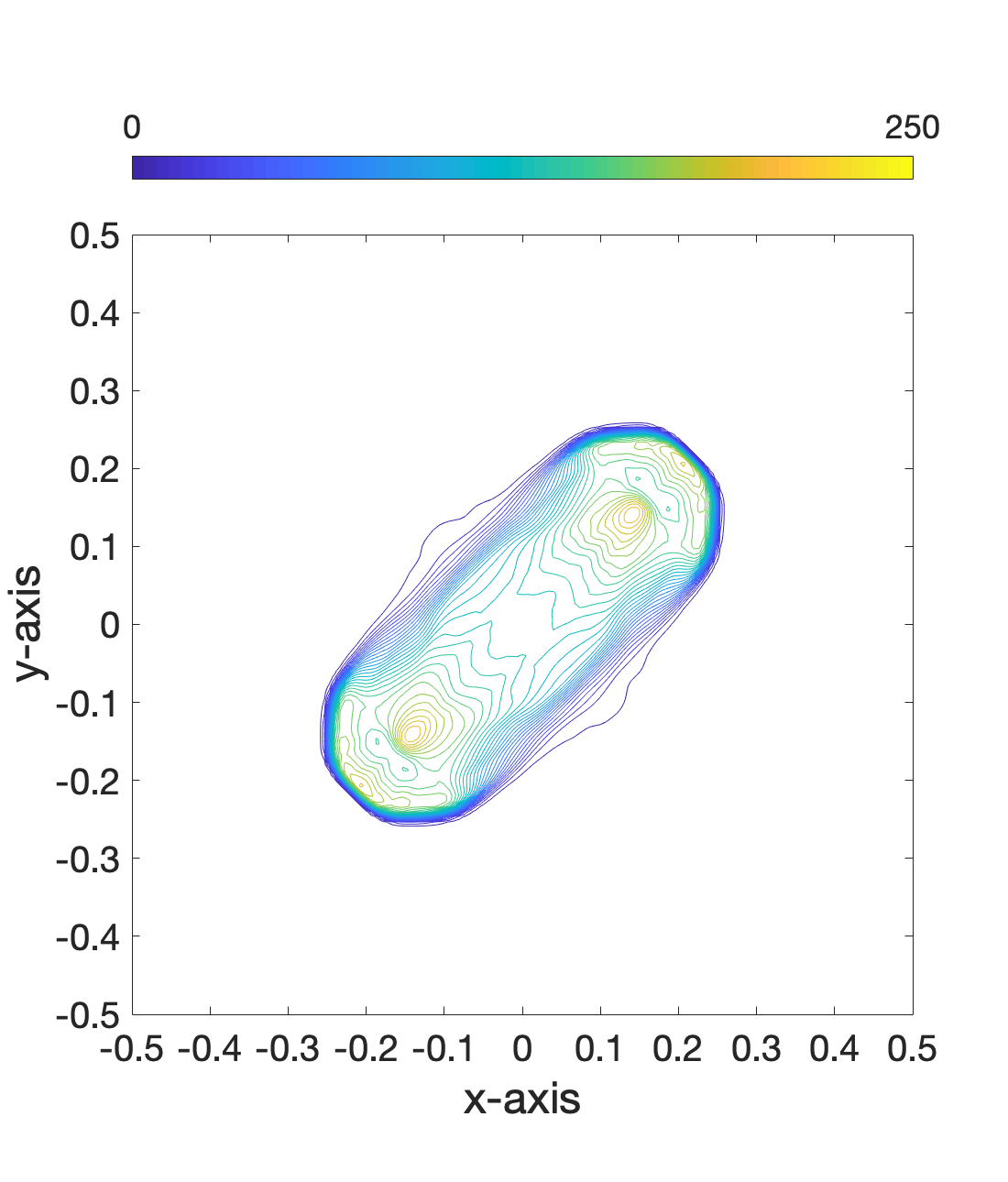}}
		\subfigure[$\| \uv \|$] {\includegraphics[width=0.45\textwidth]{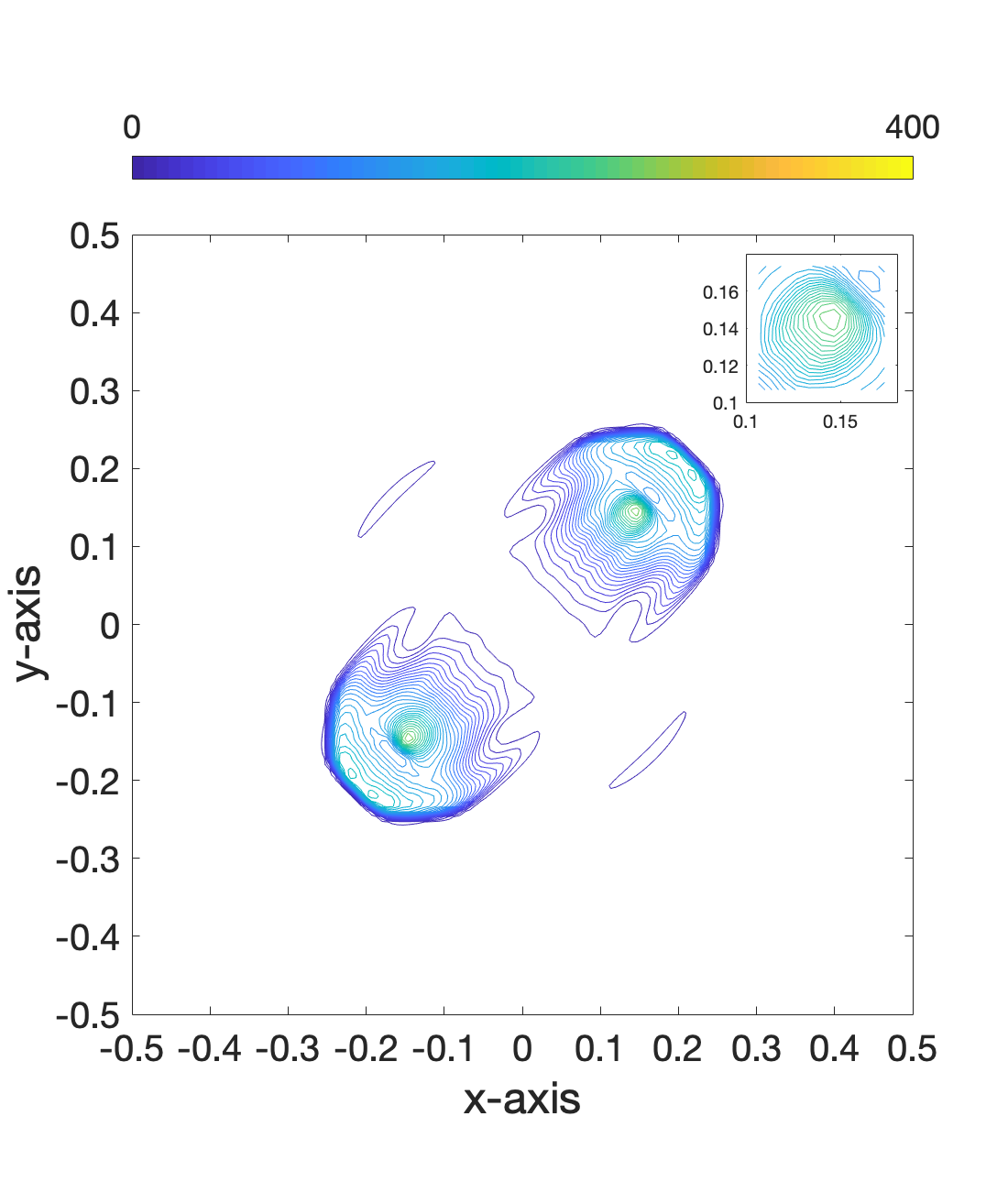}}
		\subfigure[$\| \Bv \|$] {\includegraphics[width=0.45\textwidth]{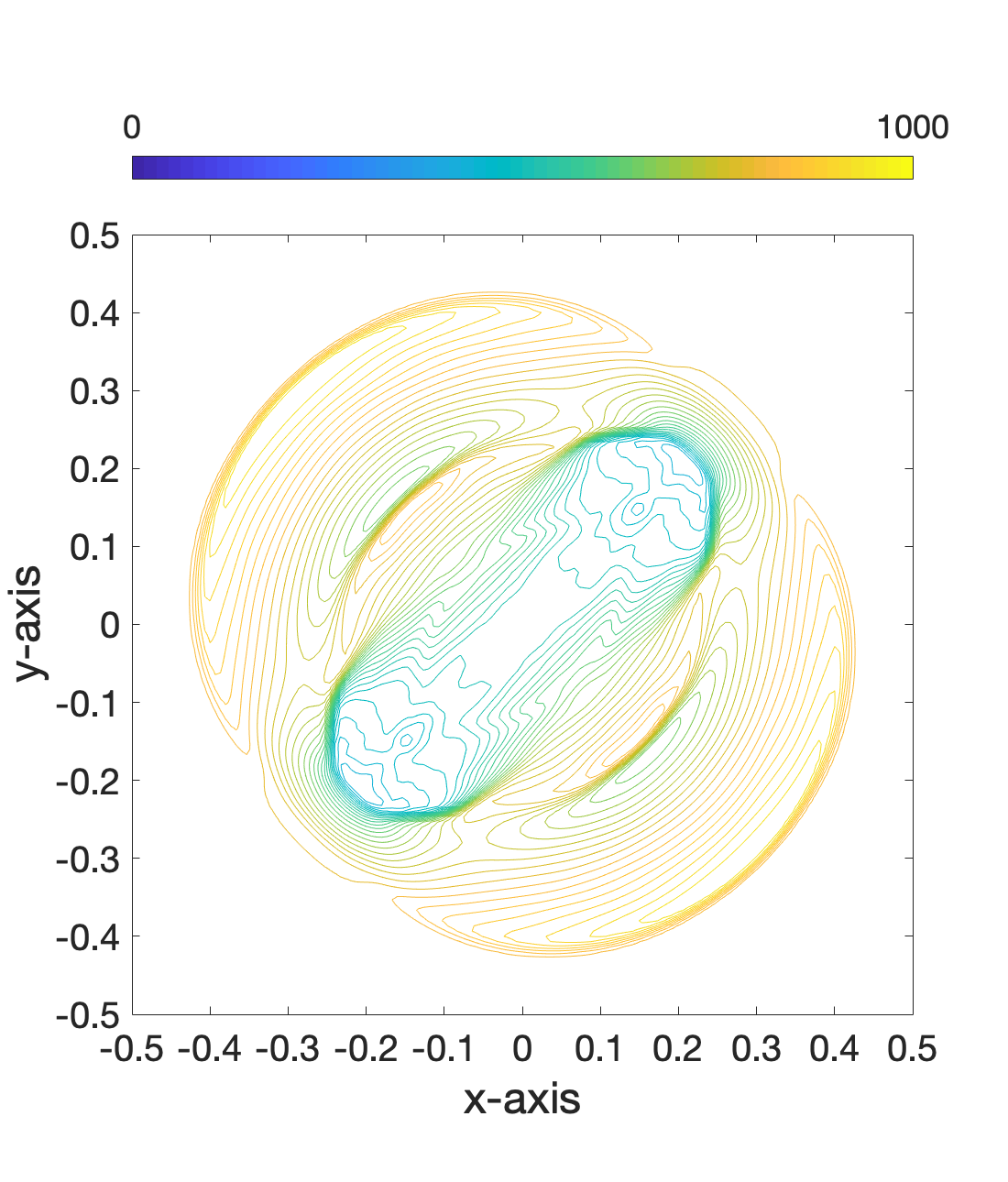}}
	\caption{\label{3dblast:diff}\em 3D Blast wave with diffusion. The max $\| \uv \|$ value is 265. }
	\end{center}
\end{figure}

In this section we demonstrate the advantage of the new approach in the context of strong shocks.  In the previous approach \cite{christlieb2014finite}, we needed to add in a diffusion limiter to stabilize the update for $\Av$ in the vicinity of a sharp change in the solution.  This had no impact on smooth solutions, and minimal impact on problems like the cloud shock.  But in the context of the blast wave problem, while the diffusion limiter stabilized the solution, comparing $\| \uv \|$ for both schemes, we have the max value of  the previous method is 261, while it is 320 for the kernel-based method. This involved decreasing the maximum values, making the solution more isotropic, and changing the structure of the contours away from the blast. In Figure \ref{3dblast:diff}, we show the solution generated by  the new method with the diffusion limiter used in our previous code. In this case, the max value of $\| \uv \|$ is 265. These results confirm that the limiter is what is causing a change in the solution structure for problems with strong shocks. Adding the diffusion limiter caused the solution to revert from the results of the new method to those obtained by the old method.
\section{Conclusion}

In this work we developed a kernel based constrained transpose scheme based on the magnetic vector potential equations in 2D and 3D  for ensuring that the solution of the magnetic filed in ideal MHD is divergence free.  The development of the method relies on a kernel based formulation of the spatial derivatives. The framework of the current method is derived from the method of lines transpose methodology and the key idea of successive convolution.  The method relies on the idea of replacing a local operator with a global operator that is as efficient as an explicit method, but because it is global it provides an unconditional stable method when coupled with explicit time stepping.  For time integration, we coupled the method with the high order explicit strong stability preserving Runge-Kutta scheme. The most important conclusion of this work is that the newly proposed method offers an approach to constrained transport that is mesh aligned for AMR and does not rely on diffusion limiter for stability in 3D, which we needed in our previous work \cite{christlieb2014finite}. Eliminating the need of diffusion limiter improves solutions where there are strong shocks, and as demonstrated in numerical simulations, the approach also works well for smooth problems.  In addition, the new method is unconditionally stable and achieves high order accuracy.  The method is robust and has been tested on a range of 2D and 3D test problems, such as the field loop and blast wave problems.


\bibliographystyle{abbrv}
\bibliography{ref}

\end{document}